\documentclass[12pt, reqno]{amsart}
\usepackage{amsmath, amstext, amsbsy, amssymb}

\newtheorem{theorem}{Theorem}[section]
\newtheorem{lemma}[theorem]{Lemma}

\newtheorem{corollary}[theorem]{Corollary}

\theoremstyle{definition}
\newtheorem{definition}[theorem]{Definition}

\theoremstyle{remark}
\newtheorem{remark}[theorem]{Remark}

\numberwithin{equation}{section}

\newcommand{\ch}{\mbox{ch} }
\newcommand{\C}{ \Bbb C }
\newcommand{\D}{ \mathcal W }
\newcommand{\DA}{{\D}(A)}
\newcommand{\End}{{\rm End}}
\newcommand{\fock}{{\mathbb H}_X}

\newcommand{\hDA}{\widehat{\D}(A)}

\newcommand{\KX}{K}
\newcommand{\lambsq}{s(\lambda)}
\newcommand{\pa}{\partial}
\newcommand{\ta}{\tau_*\alpha}
\newcommand{\tea}{\tau_*(e\alpha)}
\newcommand{\tb}{\tau_*\beta}
\newcommand{\teb}{\tau_*(e\beta)}
\newcommand{\tab}{\tau_*(\alpha\beta)}
\newcommand{\teab}{\tau_*(e\alpha\beta)}
\newcommand{\jj}{\mathfrak J}

\newcommand{\Hn}{H^*(\Xn)}
\newcommand{\Tr}{ {\rm Tr} }
\newcommand{\vac}{|0\rangle}
\newcommand{\W}{{\mathcal W}}
\newcommand{\Wax}{{\mathcal W}_X}
\newcommand{\Wbx}{{\W}^{\mathfrak B}_X}
\newcommand{\Xn}{ X^{[n]}}
\newcommand{\Z}{ \mathbb Z }

\newenvironment{demo}[1]%
{\vskip-\lastskip\medskip
  \noindent
  {\em #1.}\enspace
  }%
{\qed\par\medskip
  }

\begin{document}

\title[Hilbert schemes and $\mathcal W$ algebras]
    {Hilbert schemes and $\mathcal W$ algebras}
\author[Wei-Ping Li]{Wei-Ping Li$^1$}
\address{Department of Mathematics, HKUST, Clear Water Bay, Kowloon, Hong
Kong } \email{mawpli@ust.hk}
\thanks{${}^1$Partially supported by the grant HKUST6170/99P}

\author[Zhenbo Qin]{Zhenbo Qin$^2$}
\address{Department of Mathematics, University of Missouri, Columbia, MO
65211, USA} \email{zq@math.missouri.edu}
\thanks{${}^2$Partially supported by an NSF grant}

\author[Weiqiang Wang]{Weiqiang Wang$^3$}
\address{Department of Mathematics, University of Virginia,
Charlottesville, VA 22904} \email{ww9c@virginia.edu}
\thanks{${}^3$Partially supported by an NSF grant}


\begin{abstract}
We construct geometrically the generating fields of
a $\W$ algebra which acts irreducibly
on the direct sum of the cohomology rings of the Hilbert schemes
$\Xn$ of $n$ points on a projective surface $X$ for all $n \ge 0$.
We compute explicitly the commutators among the Fourier components of
the generating fields of the $\W$ algebra, and identify this algebra with a
$\W_{1+\infty}$-type algebra. A precise formula of certain Chern
character operators, which is essential for the construction of
the $\W$ algebra, is established in terms of the Heisenberg
algebra generators. In addition, these Chern character operators
are proved to be the zero-modes of vertex operators.
\end{abstract}

\maketitle

\date{}
\section{\bf Introduction} \label{sect_intr}

Recently, a new approach has emerged in studying the geometry of
the Hilbert schemes $\Xn$ of $n$ points on a projective surface
$X$. The starting point of this new approach is a geometric
construction, due to Nakajima \cite{Na1, Na2} and Grojnowski
\cite{Gro}, of a Heisenberg algebra which acts irreducibly on the
direct sum $\fock =\oplus_{n \ge 0} H^*(\Xn)$ of the rational
cohomology rings of $\Xn$ for all $n \ge 0$ (this was conjectured
by Vafa and Witten \cite{VW} based on G\"ottsche's formula
\cite{Got} for the Betti numbers of the Hilbert schemes $\Xn$). A
geometric construction of the Virasoro algebra acting on $\fock$
is subsequently given by Lehn \cite{Leh}. Intimate relations
between the geometry of the Hilbert schemes and vertex operators
of higher conformal weights uncovered by the authors \cite{LQW1}
provide further strong evidence toward deeper connections between
these Hilbert schemes and vertex algebras. This approach has
proved to be very fruitful in the study of the cohomology ring
structures of the Hilbert schemes $\Xn$ \cite{Leh, LQW1, LQW2,
LS1, LS2, LQW3}.

On the other hand, a distinguished class of vertex algebras, which
are called $\mathcal W$ algebras and higher-spin generalizations
of the Virasoro algebra, often appears in connection with
conformal field theories and with representation theory of affine
Kac-Moody Lie algebras (cf. the book of E.Frenkel and Ben-Zvi
\cite{FB}). Among the $\mathcal W$ algebras, a well-known example
is given by the so-called $\mathcal W_{1+\infty}$ algebra which is
the central extension of the Lie algebra of differential operators
on the circle (cf. \cite{FKRW, Kac} and the references therein).
We remark that the $\W_{1+\infty}$ algebra has an unusual feature
which is in general not true for other $\W$ algebras: the Fourier
components of the generating fields of the $\W_{1+\infty}$ algebra
are closed under the Lie bracket.

One main goal of this paper is to construct geometrically 
the generating fields of a
$\mathcal W$ (super)algebra, denoted by $\Wax$, which depends on
the projective surface $X$ and acts irreducibly on $\fock$. We
further identify the algebra $\Wax$ as an analog of the $\mathcal
W_{1+\infty}$ algebra in the framework of Hilbert schemes, which
roughly speaking, is the $\mathcal W_{1+\infty}$ algebra
parametrized by the cohomology ring $H^*(X)$ \footnote{
Just like the usual $\mathcal W_{1+\infty}$ algebra, the $\Wax$ algebra
has a dual nature: one may talk about a Lie algebra $\Wax$,
or one may talk about the vertex algebra associated to
a vacuum module of the Lie algebra $\Wax$, which is
conventionally referred to as the $\Wax$ algebra as well. The
Fourier components of the generating fields of the
vertex algebra $\Wax$ give us a linear basis for the Lie algebra
$\Wax$.}. To that end, we will
first obtain a partial description, which is of independent
interest, of the cohomology ring of the Hilbert schemes $\Xn$ for
a general projective surface $X$. It would be also interesting to
see whether or not the appearance of the $\mathcal W$ algebras in
the framework of Hilbert schemes affords an explanation from
string theory such as $S$-duality etc.

Let us explain in more details. For fixed $k\ge 0$ and $\alpha \in
H^*(X)$, we introduced in \cite{LQW1} certain cohomology class
$G_k(\alpha,n) \in H^*(\Xn)$, and then defined the Chern character
operator $\mathfrak G_k(\alpha) \in \End (\fock)$ which acts on
$H^*(\Xn)$ by the cup product with $G_k(\alpha,n)$ for each $n$
(also cf. \cite{Leh}). Such an operator approach has been
essential along this research direction. It was proved in
\cite{LQW1} that the classes $G_k(\alpha,n)$, where $0\le k<n$ and
$\alpha$ runs over a linear basis of $H^*(X)$, form a set of ring
generators for the cohomology ring $H^*(\Xn)$. Our first main
result in the present paper is to give an explicit formula for the
operator $\mathfrak G_k(\alpha)$, at least for those $\alpha$
orthogonal to the canonical class $K$, in terms of the Heisenberg
algebra generators. We emphasize that an explicit formula of this
sort contains very strong information and is technically difficult
to establish. For example, a partial information on the `leading
term' of the operator $\mathfrak G_k(\alpha)$ among other results
was essentially responsible for establishing the ring generators
statement in \cite{LQW1} mentioned above. In particular, for a
surface $X$ with numerically trivial canonical class, our precise
formula for $\mathfrak G_k(\alpha)$ provides a complete
description of the cohomology ring structure of $\Xn$ which is
totally different from the one given by Lehn and Sorger
\cite{LS2}. A corollary of our results establishes a conjecture in
\cite{LS2} (which needs to be mildly modified) concerning the
cohomology classes $G_k(\alpha,n)$. Another consequence is an
explicit formula for certain intersection numbers which were shown
earlier in \cite{LQW3} to be independent of the surface $X$.

Our construction of the algebra $\Wax$ uses the commutators of the
operators $\mathfrak G_k(\alpha)$ and the Heisenberg algebra
generators. The $\mathcal W$ algebra $\Wax$ contains both the
Heisenberg algebra of Nakajima-Grojnowski and the Virasoro algebra
of Lehn \cite{Leh} as subalgebras. We observe that a linear basis
of $\Wax$ comes from the Fourier components of explicit vertex
operators constructed from the Heisenberg vertex operators and,
most remarkably, that the operators $\mathfrak G_k(\alpha)$ for
$\alpha$ orthogonal to $K$ are precisely the zero-modes of these
vertex operators. We compute explicitly the commutation relation
among these basis elements, which in particular implies that the
algebra $\Wax$ is miraculously
closed under the Lie bracket. Note that the
commutation relations of $\Wax$ satisfy the transfer property as
formulated in \cite{LQW1}. The computation uses the explicit
formula of the operators $\mathfrak G_k(\alpha)$ obtained above
and a lengthy calculation by means of the operator product
expansion method in the theory of vertex algebras (cf. \cite{Bor,
FB, Kac}) which incorporates the transfer property.

The commutators of two basis elements in the algebra $\Wax$
typically give rise to two terms, a leading term and another term
which involves the Euler class $e$ of the projective surface $X$.
There are no third terms due to the fact that $e^2=0$, and
essentially for the same reason, central extension terms make
appearances in the commutators of $\Wax$ only for the Fourier
components of vertex operators of small conformal weights. We
further identify these leading terms as the commutators of
differential operators on the circle, and thus are justified to
regard the algebra $\Wax$ as certain topological deformation of
the $\W_{1+\infty}$ algebra in the framework of Hilbert schemes.
In particular, for a projective surface $X$ with trivial Euler
class, we have a complete identification of elements $\Wax$ with
differential operators.

Just as the results in \cite{LQW3}, our construction of the
$\mathcal W$ algebra admits a counterpart in term of the orbifold
cohomology rings of the symmetric products of a manifold, which
is worked out in \cite{QW}. In particular, when the manifold is a
point, such a counterpart specializes (with new proofs) to the
construction of the $\mathcal W_{1+\infty}$ algebra of Lascoux and
Thibon \cite{LT} using the class functions of the symmetric
groups, which in turn is an extension of a construction due to
I.~Frenkel and the third author \cite{FW} of the Virasoro algebra.

The layout of the paper is as follows. In Sect.~\ref{sect_walg},
starting from a commutative ring, we introduce a $\mathcal W$
algebra which is an analog of the $\mathcal W_{1+\infty}$ algebra.
In Sect.~\ref{sect_general}, we quickly review some known results
and constructions on the Hilbert schemes $\Xn$. In
Sect.~\ref{sect_character}, we establish the explicit formula for
the Chern character operator $\mathfrak G_k(\alpha)$ in terms of
Heisenberg generators. In Sect.~\ref{sect_walghilb}, we formulate
the main theorems on connections between $\mathcal W$ algebras and
Hilbert schemes, which are proved subsequently in
Sect.~\ref{sect_proof} by using the operator product expansion
technique.

{\bf Acknowledgment.} We thank the referee for a careful
reading and pointing out an error in the original proof
of Theorem 4.7.
\section{\bf The definition of $\W$ algebras} \label{sect_walg}

Let $A =\oplus_{n=0}^d A_n$ (where $A_d \cong \C$) be a
finite-dimensional graded complex vector space. We write $|\alpha|
=n$ if $\alpha \in A_n$. We define a super structure (i.e. the
$\Z/2\Z$-grading) on $A$ by letting $A_{\bar{0}} =\oplus_{n\in
2\Z_+} A_n$ and $A_{\bar{1}} =\oplus_{n\in 1+ 2\Z_+} A_n$ (here
$\Z_+$ stands for the set of all nonnegative integers). We assume
that $A$ affords an algebra structure which is super commutative
and compatible with the $\Z_+$-grading. We further assume that
there exists a linear operator $\Tr :A \rightarrow \C$ which is
zero on the graded subspace $A_n$ unless $n$ is the top degree.
This defines a supercommutative bilinear form $(-,-):A\times A
\rightarrow \C$ by $(\alpha,\beta) =\Tr(\alpha \beta)$.

Let $t$ be an indeterminate and let
$\partial_t=\displaystyle{\frac{d}{dt}}$. Let ${\D}_{as}$ be the
associative algebra of regular differential operators on the
circle $S^1$. Denote by ${\DA}_{as}$ the associative superalgebra
${\D}_{as} \otimes A$. It has a linear basis given by
\begin{eqnarray*}
  J^p_k (\alpha)= - t^{p+k} ( \partial_t )^p \otimes \alpha, \quad
    p \in \Z_{+}, k \in \Z, \alpha \in A.
\end{eqnarray*}
A different linear basis of $\DA_{as}$ is given by
\begin{eqnarray*}
 L^p_k (\alpha) =
 - t^{k} D^p \otimes \alpha,  \quad p \in \Z_{+}, k \in \Z, \alpha \in A,
\end{eqnarray*}
where $D = t \partial_t$. Note that $f(D) t = t f( D +1)$ for
every polynomial $f(w) \in \C [w]$. Hence, we have $J^p_k(\alpha)
= -t^k D(D-1)\cdots (D-p+1) \otimes \alpha$.

Let $\DA$ denote the Lie superalgebra obtained from the
associative superalgebra ${\DA}_{as}$ by taking the usual super
bracket of operators:
\begin{eqnarray*}
[X\otimes \alpha, Y \otimes \beta] = (XY -YX) \otimes (\alpha\beta
).
\end{eqnarray*}
Assuming that $\alpha$ and $\beta$ are of homogenous degree, we
see that $[X\otimes \alpha, Y \otimes \beta]
=(-1)^{1+|\alpha|\cdot | \beta|}[Y \otimes \beta, X\otimes
\alpha]$. The commutation relation in $\DA$ is given by
\begin{eqnarray}
 && \left[
     t^r f(D) \otimes \alpha, t^s g(D) \otimes \beta
  \right] \nonumber \\
    & = & t^{r+s}
    \left(
      f(D + s)g(D) - f(D)g(D +r)
    \right) \otimes (\alpha\beta) .
  \label{eq_com}
\end{eqnarray}

When $A =\C$ (hence $d=0$), we will simply write the Lie
superalgebra $\DA$ as $\D$, which is the usual Lie algebra of
differential operator on the circle, cf. e.g. \cite{FKRW, Kac}.
The algebra $\D$ affords a universal central extension which is
usually referred to as the $\W_{1+\infty}$ algebra. In
$\W_{1+\infty}$, the operators $L^0_k$ generate a Heisenberg
algebra, while the operators $L^1_k$ generate a Virasoro algebra.

The algebra $\DA$ has a natural {\em weight} filtration
\begin{eqnarray*}
\DA^0 \subset\DA^1 \subset \DA^2 \subset \ldots \subset \DA
\end{eqnarray*}
by letting the weight of $L^p_k(\alpha)$ be $p$. Clearly $[\DA^p,
\DA^q]\subset \DA^{p+q-1}$. We will denote the associated graded
algebra by $\mathcal{GW}(A)$. The leading term (according to the
weight filtration) of $J^p_k(\alpha)$ is just $L^p_k(\alpha)$, and
so $J^p_k(\alpha)$ and $L^p_k(\alpha)$ give rise to the same
element in the graded algebra $\mathcal{GW}(A)$ which is denoted
by $\mathfrak L^p_k(\alpha)$. We easily derive the following
commutation relation from (\ref{eq_com}):
\begin{eqnarray*}
  [\mathfrak L^p_m(\alpha), \mathfrak L^q_n(\beta)]
   & = &      (q m -p n) \cdot \mathfrak L^{p+q-1}_{m+n}(\alpha\beta).
\end{eqnarray*}

The $\W$ (super)algebra $\hDA$ is a central extension of the Lie
superalgebra $\mathcal{GW}(A)$ by a one-dimensional center with a
specified generator $C$:
\begin{eqnarray*}
 0 \longrightarrow \C C \longrightarrow \hDA \longrightarrow
 \mathcal{GW}(A) \longrightarrow 0,
\end{eqnarray*}
such that the commutators in $\hDA =\mathcal{GW}(A)+\C C$ are
given by:
\begin{eqnarray}  \label{eq_gradedcomm}
  [\mathfrak L^p_m(\alpha), \mathfrak L^q_n(\beta)]
   & = &
 \left\{
     \everymath{\displaystyle}
      \begin{array}{ll}
        m\delta_{m,-n} \, \Tr(\alpha\beta) \cdot C, &{\rm if} \;p=q=0, \\
       (q m -p n) \cdot \mathfrak L^{p+q-1}_{m+n}(\alpha\beta), &{\rm otherwise}.
      \end{array}
 \right.
\end{eqnarray}
In the above, we have used $\delta_{a,b}$ to denote $1$ if $a=b$
and $0$ otherwise.

In this paper, we will be mainly interested in the $\W$
(super)algebra $\hDA$ when $A$ is (a subring of) the cohomology
ring $H^*(X)$ of a projective surface $X$, with the trace being
defined by $\Tr(\alpha) =-\displaystyle{\int_X\alpha}$ for $\alpha
\in H^*(X)$. Here and below $H^*(X)$ always denote the cohomology
group/ring of $X$ with rational coefficient.
\section{\bf Basics on Hilbert schemes of points on surfaces} \label{sect_general}

Let $X$ be a smooth projective complex surface with the canonical
class $\KX$ and the Euler class $e$, and $\Xn$ be the Hilbert
scheme of points in $X$. An element in $\Xn$ is represented by a
length-$n$ $0$-dimensional closed subscheme $\xi$ of $X$. For $\xi
\in \Xn$, let $I_{\xi}$ be the corresponding sheaf of ideals. It
is well known that $\Xn$ is smooth. Sending an element in $\Xn$ to
its support in the symmetric product ${\rm Sym}^n(X)$, we obtain
the Hilbert-Chow morphism $\pi_n: \Xn \rightarrow {\rm Sym}^n(X)$,
which is a resolution of singularities. Define the universal
codimension-$2$ subscheme:
\begin{eqnarray*}
{ \mathcal Z}_n=\{(\xi, x) \subset \Xn\times X \, | \, x\in
 {\rm Supp}{(\xi)}\}\subset \Xn\times X.
\end{eqnarray*}
Denote by $p_1$ and $p_2$ the projections of $\Xn \times X$ to
$\Xn$ and $X$ respectively. Let
\begin{eqnarray*}
\fock = \oplus_{n=0}^\infty \Hn
\end{eqnarray*}
be the direct sum of total cohomology groups of the Hilbert
schemes $\Xn$.

For $m \ge 0$ and $n > 0$, let $Q^{[m,m]} = \emptyset$ and define
$Q^{[m+n,m]}$ to be the closed subset:
$$\{ (\xi, x, \eta) \in X^{[m+n]} \times X \times X^{[m]} \, | \,
\xi \supset \eta \text{ and } \mbox{Supp}(I_\eta/I_\xi) = \{ x \}
\}.$$

We recall Nakajima's definition of the Heisenberg operators
\cite{Na1}. Let $n > 0$. The linear operator $\mathfrak
a_{-n}(\alpha) \in \End(\fock)$ with $\alpha \in H^*(X)$ is
defined by
$$\mathfrak a_{-n}(\alpha)(a) = \tilde{p}_{1*}([Q^{[m+n,m]}] \cdot
\tilde{\rho}^*\alpha \cdot \tilde{p}_2^*a)$$
for $a \in H^*(X^{[m]})$, where $\tilde{p}_1, \tilde{\rho},
\tilde{p}_2$ are the projections of $X^{[m+n]} \times X \times
X^{[m]}$ to $X^{[m+n]}, X, X^{[m]}$ respectively. Define
$\mathfrak a_{n}(\alpha) \in \End(\fock)$ to be $(-1)^n$ times the
operator obtained from the definition of $\mathfrak
a_{-n}(\alpha)$ by switching the roles of $\tilde{p}_1$ and $
\tilde{p}_2$. We often refer to $\mathfrak a_{-n}(\alpha)$ (resp.
$\mathfrak a_n(\alpha)$) as the {\em creation} (resp. {\em
annihilation})~operator. We also set $\mathfrak a_0(\alpha) =0$.

For $n > 0$ and a homogeneous class $\gamma \in H^*(X)$, let
$|\gamma| = s$ if $\gamma \in H^s(X)$, and let $G_i(\gamma, n)$ be
the homogeneous component in $H^{|\gamma|+2i}(\Xn)$ of
\[
 G(\gamma, n) = p_{1*}(\ch({\mathcal O}_{{\mathcal Z}_n}) \cdot
 p_2^*{\rm td}(X) \cdot p_2^*\gamma) \in \Hn
\]
where $\ch({\mathcal O}_{{\mathcal Z}_n})$ denotes the Chern
character of the structure sheaf ${\mathcal O}_{{\mathcal Z}_n}$
and ${\rm td}(X) $ denotes the Todd class. Here and below we omit
the Poincar\'e duality used to switch a homology class to a
cohomology class and vice versa. We extend the notion $G_i(\gamma,
n)$ linearly to an arbitrary class $\gamma \in H^*(X)$.
We also set $G(\gamma, 0) =0$.

The {\it Chern character operator} ${\mathfrak G}_i(\gamma) \in
\End({\fock})$ is defined to be the operator acting on the
component $H^*(\Xn)$ by the cup product with $G_i(\gamma, n)$. It
was proved in \cite{LQW1} that the cohomology ring of $\Xn$ is
generated by the classes $G_{i}(\gamma, n)$ where $0 \le i < n$
and $\gamma$ runs over a linear basis of $H^*(X)$. Let $\mathfrak
d = \mathfrak G_1(1_X)$ where $1_X$ is the fundamental cohomology
class of $X$. The operator $\mathfrak d$ was first introduced in
\cite{Leh}.
For a linear operator $\mathfrak f \in \End(\fock)$, define its
{\it derivative} $\mathfrak f'$ by $\mathfrak f' = [\mathfrak d,
\mathfrak f]$. The higher derivative $\mathfrak f^{(k)}$ is
defined inductively by $\mathfrak f^{(k)} = [\mathfrak d,
\mathfrak f^{(k-1)}]$.

Let $:\frak a_{m_1}\frak a_{m_2}:$ be $\frak a_{m_1}\frak a_{m_2}$
when $m_1 \le m_2$ and $\frak a_{m_2}\frak a_{m_1}$ when $m_1 >
m_2$. For $n \in \Z$, define a linear map $\mathfrak L_n: H^*(X)
\to \End(\fock)$ by $\displaystyle{\mathfrak L_n = -\frac{1}{2}
\cdot \sum_{m \in \Z} :\frak a_{m} \frak a_{n-m}: \tau_{2*}}$.
Here for $k \ge 1$, $\tau_{k*}: H^*(X) \to H^*(X^k)$ is the linear
map induced by the diagonal embedding $\tau_k: X \to X^k$, and
$\mathfrak a_{m_1} \cdots \mathfrak a_{m_k}(\tau_{k*}(\alpha))$
denotes $\sum_j \mathfrak a_{m_1}(\alpha_{j,1}) \cdots \mathfrak
a_{m_k}(\alpha_{j,k})$
when $\tau_{k*}\alpha = \sum_j \alpha_{j,1} \otimes \cdots \otimes
\alpha_{j, k}$ via the K\"unneth decomposition of $H^*(X^k)$.

The following is a combination of various theorems from \cite{Na2,
Gro, Leh, LQW1}. Our notations and convention of signs are
consistent with \cite{LQW3}.

\begin{theorem} \label{commutator}
Let $k \ge 0, n,m \in \Z$ and $\alpha, \beta \in H^*(X)$. Then,

\begin{enumerate}
\item[{\rm (i)}] The operators $\mathfrak a_n(\alpha)$ satisfy
a Heisenberg algebra commutation relation:
\begin{eqnarray*}
\displaystyle{[\mathfrak a_m(\alpha), \mathfrak a_n(\beta)] = -m
\; \delta_{m,-n} \int_X(\alpha \beta) \cdot {\rm Id}_{\fock}}.
\end{eqnarray*}
The space $\fock$ is an irreducible module over the Heisenberg
algebra generated by the operators $\mathfrak a_n(\alpha)$ with a
highest~weight~vector $\vac=1 \in H^0(X^{[0]}) \cong \mathbb Q$.

\item[{\rm (ii)}] $[\mathfrak L_m(\alpha), \mathfrak a_n(\beta)]
= -n \cdot \mathfrak a_{m+n}(\alpha \beta)$;

\item[{\rm (iii)}] $\displaystyle{\mathfrak a_n'(\alpha)
= n \cdot \mathfrak L_n(\alpha) - \frac{n(|n|-1)}{2} \cdot
\mathfrak a_n(\KX \alpha)}$;

\item[{\rm (iv)}] The operators $\mathfrak L_n(\alpha)$ satisfy
the Virasoro algebra commutation relation:
\begin{eqnarray*}
[\mathfrak L_m(\alpha),\mathfrak L_n(\beta)]
 = (m-n) \cdot \mathfrak L_{m+n}(\alpha\beta) +\frac{m^3-m}{12}\,
 \delta_{m,-n} \int_X (e \alpha\beta) \cdot {\rm Id}_{\fock}.
\end{eqnarray*}

\item[{\rm (v)}] $\displaystyle{[\mathfrak G_k(\alpha),
\mathfrak a_{-1}(\beta)] = \frac{1}{k!} \cdot \mathfrak
a_{-1}^{(k)}(\alpha \beta)}$.
\end{enumerate}
\end{theorem}

The Lie brackets in the above theorem are understood in the super
sense according to the parity of the cohomology degrees of the
cohomology classes involved. Also, it follows from
Theorem~\ref{commutator}~(i) that the space $\fock$ is linearly
spanned by all the Heisenberg monomials $\mathfrak
a_{n_1}(\alpha_1) \cdots \mathfrak a_{n_k}(\alpha_k)\cdot \vac$
where $k \ge 0$ and $n_1, \ldots, n_k < 0$.

We will need later on the following lemma proved in \cite{LQW3}.

\begin{lemma} \label{tau_k_tau_{k-1}}
Let $k, s \ge 1$, $n_1, \ldots, n_k, m_1, \ldots, m_s \in \Z$, and
$\alpha, \beta \in H^*(X)$.
 \begin{enumerate}
\item[{\rm (i)}] The commutator
$[\mathfrak a_{n_1} \cdots \mathfrak a_{n_{k}} (\tau_{k*}\alpha),
\mathfrak a_{m_1} \cdots \mathfrak a_{m_{s}}(\tau_{s*}\beta)]$ is
equal to
\begin{eqnarray*}
-\sum_{t=1}^k \sum_{j=1}^s n_t \delta_{n_t,-m_j} \cdot \left (
\prod_{\ell=1}^{j-1} \mathfrak a_{m_\ell} \prod_{1 \le u \le k, u
\ne t} \mathfrak a_{n_u} \prod_{\ell=j+1}^{s} \mathfrak a_{m_\ell}
\right )(\tau_{(k+s-2)*}(\alpha\beta));
\end{eqnarray*}

\item[{\rm (ii)}] The derivative $(\mathfrak a_{n_1} \cdots \mathfrak
a_{n_k}(\tau_{k*}\alpha))'$ is equal to
\begin{eqnarray*}
 & &-\sum_{j=1}^k \frac{n_j}{2} \cdot \sum_{m_1 + m_2 = n_j}
    \mathfrak a_{n_1} \cdots \mathfrak a_{n_{j-1}}
    :\mathfrak a_{m_1}\mathfrak a_{m_2}:
    \mathfrak a_{n_{j+1}} \cdots \mathfrak a_{n_k}(\tau_{(k+1)*}\alpha)   \\
 & &- \sum_{j=1}^k \frac{n_j(|n_j|-1)}{2} \cdot
    \mathfrak a_{n_1} \cdots \mathfrak a_{n_k}(\tau_{k*}(\KX\alpha));
\end{eqnarray*}

\item[{\rm (iii)}] Let $j$ satisfy $1 \le j < k$. Then,
$\mathfrak a_{n_1} \cdots \mathfrak a_{n_k}(\tau_{k*}\alpha)$ is
equal to
\begin{eqnarray*}
\left ( \prod_{1 \le s < j} \mathfrak a_{n_s} \cdot \mathfrak
a_{n_{j+1}} \mathfrak a_{n_{j}} \cdot \prod_{j+1 < s \le k}
\mathfrak a_{n_s} \right ) (\tau_{k*}\alpha) - n_j
\delta_{n_j,-n_{j+1}} \prod_{1 \le s \le k \atop s \ne j, j+1}
\mathfrak a_{n_s}(\tau_{(k-2)*}(e\alpha)).
\end{eqnarray*}
 \end{enumerate}
\end{lemma}

In the lemma above and throughout the paper, 
the products of Heisenberg operators are understood
in the increasing order of the parametrizing indices
from the left to the right,
e.g., $\prod_{1 \le s < j} \mathfrak a_{n_s}
= \mathfrak a_{n_1}\mathfrak a_{n_2} \ldots 
\mathfrak a_{n_{j-1}}$ and
$\prod_{1 \le u \le k, u
\ne t} \mathfrak a_{n_u} =
\mathfrak a_{n_1} \ldots\mathfrak a_{n_{t-1}}
\mathfrak a_{n_{t+1}} \ldots\mathfrak a_{n_k}.$

Note that the commutator in Lemma~\ref{tau_k_tau_{k-1}}~(i)
satisfies the {\em transfer property}, that is, it depends only on
the cup product $(\alpha\beta)$. Such a property was first
formulated and emphasized in \cite{LQW1}, and will also be
incorporated into the formulation of the operator product
expansions used later in this paper.

\section{\bf Chern character operators}
 \label{sect_character}

In this section, we determine the operators $\mathfrak
a_n^{(k)}(\alpha)$ and $\mathfrak G_k(\alpha)$ when the cohomology
class $\alpha$ is orthogonal to the canonical class $K$. In
addition, we obtain general expressions of $\mathfrak
a_n^{(k)}(\alpha)$ and $\mathfrak G_k(\alpha)$ when $\alpha \in
H^*(X)$ is arbitrary.
\subsection{Formulas for the derivatives of Heisenberg operators}

\begin{definition} \label{partition}
Let $X$ be a smooth projective surface.

\begin{enumerate}
\item[{\rm (i)}] Let $\alpha \in H^*(X)$, and $\lambda = (\cdots
(-2)^{m_{-2}}(-1)^{m_{-1}} 1^{m_1}2^{m_2} \cdots)$ be a {\em
generalized partition} of the integer $n = \sum_i i m_i$ whose
part $i\in \Z$ has multiplicity $m_i$. Define $\ell(\lambda) =
\sum_i m_i$, $|\lambda| = \sum_i i m_i = n$, $\lambsq  = \sum_i
i^2 m_i$, $\lambda^! = \prod_i m_i!$, and
\begin{eqnarray*}
\mathfrak a_{\lambda}(\tau_*\alpha) = \left ( \prod_i (\mathfrak
a_i)^{m_i} \right ) (\tau_{\ell(\lambda)*}\alpha)
\end{eqnarray*}
where the product $\prod_i (\mathfrak
a_i)^{m_i} $ is understood to be 
$\cdots \mathfrak a_{-2}^{m_{-2}} \mathfrak a_{-1}^{m_{-1}} 
 \mathfrak a_{1}^{m_{1}} \mathfrak a_{2}^{m_{2}}\cdots.$
Let $-\lambda$ be the generalized partition whose multiplicity of
$i \in \Z$ is $m_{-i}$.

\item[{\rm (ii)}] A generalized partition becomes a {\em partition}
in the usual sense if the multiplicity $m_ i = 0$ for every $i <
0$. A partition $\lambda$ of $n$ is denoted by $\lambda \vdash n$.

\item[{\rm (iii)}] Define $\mathfrak A_X =\{\alpha\in H^*(X)|
\KX \alpha = 0\}$ which is an ideal in the ring $H^*(X)$.
\end{enumerate}
\end{definition}

Our starting point is the following theorem.

\begin{theorem} \label{derivative2.1}
Let $k \ge 0$, $n \in \Z$, and $\alpha\in \mathfrak A_X$. Then,
$\mathfrak a_n^{(k)}(\alpha)$ is equal to
\begin{eqnarray} \label{derivative2.1.1}
\quad(-n)^k k! \left ( \sum_{\ell(\lambda) = k+1, |\lambda|=n}
\frac{1}{\lambda^!} \mathfrak a_{\lambda}(\tau_{*}\alpha) -
\sum_{\ell(\lambda) = k-1, |\lambda|=n} \frac{\lambsq -1}{24
\lambda^!} \mathfrak a_{\lambda}(\tau_{*}(e\alpha)) \right ).
\end{eqnarray}
\end{theorem}
\begin{demo}{Proof}
Use induction on $k$. When $k = 0$ or $1$, the theorem is
trivially true. Next, assume that the theorem is true for
$\mathfrak a_n^{(k)}(\alpha)$, i.e., $\mathfrak a_n^{(k)}(\alpha)$
is given by (\ref{derivative2.1.1}). By
Lemma~\ref{tau_k_tau_{k-1}} (ii) and (iii), we conclude from
(\ref{derivative2.1.1}) that $\mathfrak a_n^{(k+1)}(\alpha)$ is of
the form
\begin{eqnarray*}
\sum_{\ell(\lambda) = k+2, |\lambda|=n} {\tilde f}_{1_X}(\lambda)
\mathfrak a_{\lambda}(\tau_{*}\alpha) + \sum_{\ell(\lambda) = k,
|\lambda|=n} {\tilde f}_{e}(\lambda) \mathfrak
a_{\lambda}(\tau_{*}(e\alpha)).
\end{eqnarray*}
Fix a generalized partition $\lambda = (\cdots
(-2)^{m_{-2}}(-1)^{m_{-1}} 1^{m_1}2^{m_2} \cdots)$ with
$|\lambda|=n$. In the following, we compute the coefficients
${\tilde f}_{1_X}(\lambda)$ and ${\tilde f}_{e}(\lambda)$
separately.

(i) We start with the computation of ${\tilde f}_{1_X}(\lambda)$.
Note that $\ell(\lambda) = k+2$. Going from $\mathfrak
a_n^{(k)}(\alpha)$ to $\mathfrak a_n^{(k+1)}(\alpha) = (\mathfrak
a_n^{(k)}(\alpha))'$, we see that those terms in $\mathfrak
a_n^{(k)}(\alpha)$ whose derivatives contain the term $\mathfrak
a_{\lambda}(\tau_{*}\alpha)$ in $\mathfrak a_n^{(k+1)}(\alpha)$
are of the form $ \mathfrak a_{\lambda_{i, j}}(\tau_{*}\alpha)$
where $i \le j$, $m_i \ge 1$, $m_j \ge 1$, and $\lambda_{i, j}$
stands for the generalized partition obtained from $\lambda$ by
subtracting $1$ from the multiplicities of $i$ and $j$ and by
adding $1$ to the multiplicity of $(i+j)$. For fixed $i$ and $j$,
$\{ \mathfrak a_{\lambda_{i, j}}(\tau_{*}\alpha) \}'$ is equal to
\begin{eqnarray*}
& &\left (- \frac{i+j}{2} \right ) (m_{i+j}+1) (2 - \delta_{i, j})
   \mathfrak a_{\lambda}(\tau_{*}\alpha)    \\
&+&\sum_{\mu \ne \lambda \atop \ell(\mu) = k+2, |\mu|=n}
   d_{\mu} \mathfrak a_{\mu}(\tau_{*}\alpha) +
\sum_{\ell(\mu) = k, |\mu|=n}
   d_{\mu} \mathfrak a_{\mu}(\tau_{*}(e\alpha))
\end{eqnarray*}
for some constant $d_{\mu}$.
By (\ref{derivative2.1.1}), the coefficient of $\mathfrak
a_{\lambda_{i, j}}(\tau_{*}\alpha)$ in $\mathfrak
a_n^{(k)}(\alpha)$ is
\begin{eqnarray*}
{(-n)^k k! \over (\lambda_{i, j})^!} = {(-n)^k k! \over \lambda^!}
\cdot {m_i(m_j-\delta_{i, j}) \over m_{i+j}+1}.
\end{eqnarray*}

Combining all of these, we see that the coefficient of $\mathfrak
a_{\lambda}(\tau_{*}\alpha)$ in $\mathfrak a_n^{(k+1)}(\alpha)$ is
\begin{eqnarray*}
& &\sum_{i \le j} {(-n)^k k! \over \lambda^!} \cdot
{m_i(m_j-\delta_{i, j})
   \over m_{i+j}+1}
   \left (- {i+j \over 2} \right ) (m_{i+j}+1) (2 - \delta_{i, j})  \\
&=&{(-n)^k k! \over \lambda^!} \left (
   -\sum_{i < j} m_im_j (i+j) - \sum_{i} m_i(m_i-1)i \right )  \\
&=&{(-n)^k k! \over \lambda^!} \left (
   -{1 \over 2} \sum_{i, j} m_im_j (i+j)
   + \sum_{i} m_i i \right )
={(-n)^{k+1} (k+1)! \over \lambda^!}
\end{eqnarray*}
where we have used $\sum_i m_i = \ell(\lambda) = k+2$ and $\sum_i
im_i = |\lambda| = n$ in the last step.

(ii) Now we compute ${\tilde f}_{e}(\lambda)$. In this case,
$\ell(\lambda) = k$. We shall prove that the coefficient of
$\mathfrak a_{\lambda}(\tau_{*}(e\alpha))$ in $\mathfrak
a_n^{(k+1)}(\alpha)$ is equal to $\displaystyle{{(-n)^{k+1} (k+1)!
\over \lambda^!} c}$ where
\begin{eqnarray*}
c = {1 - \lambsq \over 24} = {1 - \sum_i i^2 m_i \over 24}.
\end{eqnarray*}
In view of Lemma~\ref{tau_k_tau_{k-1}} (ii) and (iii), there are
exactly two sources contributing to the term $\mathfrak
a_{\lambda}(\tau_{*}(e\alpha))$ in $\mathfrak a_n^{(k+1)}(\alpha)$
from the derivatives of the terms in $\mathfrak
a_n^{(k)}(\alpha)$. In the following, we handle these two
different sources separately.

The first is similar to (i) above, and comes from the derivatives
$\{ \mathfrak a_{\lambda_{i, j}}(\tau_{*}(e\alpha)) \}'$ where $i
\le j$, $m_i \ge 1$, $m_j \ge 1$. For fixed $i$ and $j$, $\{
\mathfrak a_{\lambda_{i, j}}(\tau_{*}(e\alpha)) \}'$ is equal to
\begin{eqnarray*}
\left (- {i+j \over 2} \right ) (m_{i+j}+1) (2 - \delta_{i, j})
\mathfrak a_{\lambda}(\tau_{*}(e\alpha)) + \sum_{\mu \ne \lambda
\atop \ell(\mu) = k, |\mu|=n} d_{\mu} \mathfrak
a_{\mu}(\tau_{*}(e\alpha))
\end{eqnarray*}
for some constant $d_{\mu}$.
By (\ref{derivative2.1.1}), the coefficient of $\mathfrak
a_{\lambda_{i, j}}(\tau_{*}(e\alpha))$ in $\mathfrak
a_n^{(k)}(\alpha)$ is equal to
\begin{eqnarray*}
{(-n)^k k! \over (\lambda_{i, j})^!} \cdot {1 - s(\lambda_{i, j})
\over 24} = {(-n)^k k! \over \lambda^!} \cdot {m_i(m_j-\delta_{i,
j}) \over m_{i+j}+1} \cdot \left (c - {ij \over 12} \right ).
\end{eqnarray*}
Combining all of these and noting $\sum_i m_i = k$ and $\sum_i
im_i= n$, we see that the total contribution of the first source
to the coefficient of $\mathfrak a_{\lambda} (\tau_{*}(e\alpha))$
in $\mathfrak a_n^{(k+1)}(\alpha)$ is
\begin{eqnarray} \label{derivative2.2}
& &\sum_{i \le j} {(-n)^k k! \over \lambda^!} \cdot
{m_i(m_j-\delta_{i, j})
   \over m_{i+j}+1} \cdot \left (c - {ij \over 12} \right )
   \left (- {i+j \over 2} \right ) (m_{i+j}+1) (2 - \delta_{i, j})
   \nonumber \\
&=&{(-n)^k k! \over \lambda^!} \left ( nc-knc
   + {n \over 12} \sum_i i^2m_i -{1 \over 12} \sum_i i^3m_i \right ).
\end{eqnarray}

The second source comes from the derivatives $\{ \mathfrak
a_{\tilde \lambda_{i, j}}(\tau_{*}\alpha) \}'$ where $i \ne 0$, $0
< 2j \le |i|$, $m_i \ge 1$, and $\tilde \lambda_{i, j}$ is defined
as follows. There are two cases depending on $i>0$ or $i<0$. For
$i>0$, define $\tilde \lambda_{i, j}$ to be the generalized
partition obtained from $\lambda$ by subtracting $1$ from the
multiplicity of $i$ and by adding $1$ to the multiplicities of $j$
and $(i-j)$. Then, $\{ \mathfrak a_{\tilde \lambda_{i,
j}}(\tau_{*}\alpha) \}'$ contains a term $(\ldots \mathfrak
a_j\ldots \mathfrak a_{i-j}'\ldots)$ which contains a term
$(\ldots \mathfrak a_j \ldots \mathfrak a_{-j} \mathfrak
a_i\ldots)$ after expanding the derivative $\mathfrak a_{i-j}'$.
When we switch $\mathfrak a_j$ with $\mathfrak a_{-j}$, we get the
term $\mathfrak a_{\lambda}(\tau_{*}(e\alpha))$ in $\mathfrak
a_n^{(k+1)} (\alpha)$. For $i<0$, define $\tilde \lambda_{i, j}$
to be the generalized partition obtained from $\lambda$ by
subtracting $1$ from the multiplicity of $i$ and by adding $1$ to
the multiplicities of $-j$ and $(i+j)$. Then $\{ \mathfrak
a_{\tilde \lambda_{i, j}}(\tau_{*}\alpha) \}'$ contains a term
$(\ldots \mathfrak a_{i+j}'\ldots \mathfrak a_{-j}\ldots)$ which
contains a term $(\ldots \mathfrak a_i \mathfrak a_{j}\ldots
\mathfrak a_{-j}\ldots)$ after expanding $\mathfrak a_{i+j}'$.
Switching $\mathfrak a_j$ with $\mathfrak a_{-j}$, we get
$\mathfrak a_{\lambda}(\tau_{*}(e\alpha))$ in $\mathfrak
a_n^{(k+1)} (\alpha)$.

More precisely, for fixed $i > 0$ and $j$, the derivative $\{
\mathfrak a_{\tilde \lambda_{i, j}}(\tau_{*}\alpha) \}'$ is
\begin{eqnarray*}
& &j(i-j) \cdot {(m_{i-j}+1)(m_j+1+\delta_{i, 2j}) \over 1+
\delta_{i, 2j}}
   \cdot \mathfrak a_{\lambda}(\tau_{*}(e\alpha)) +  \\
&+&\sum_{\mu \ne \lambda \atop \ell(\mu) = k, |\mu|=n}
   d_{\mu} \mathfrak a_{\mu}(\tau_{*}(e\alpha)) +
   \sum_{\ell(\mu) = k+2, |\mu|=n} d_{\mu} \mathfrak a_{\mu}(\tau_{*}\alpha)
\end{eqnarray*}
for some constant $d_{\mu}$.
By (\ref{derivative2.1.1}), the coefficient of $\mathfrak
a_{\tilde \lambda_{i, j}}(\tau_{*}(e\alpha))$ in $\mathfrak
a_n^{(k)}(\alpha)$ is equal to
\begin{eqnarray*}
{(-n)^k k! \over (\tilde \lambda_{i, j})^!} = {(-n)^k k! \over
\lambda^!} \cdot {m_i \over (m_{i-j}+1)(m_j+1+\delta_{i, 2j})}.
\end{eqnarray*}
It follows that the total contribution of the second source with
$i > 0$ to the coefficient of the term $\mathfrak
a_{\lambda}(\tau_{*}(e\alpha))$ in $\mathfrak a_n^{(k+1)}(\alpha)$
is equal to
\begin{eqnarray} \label{derivative2.3}
& &\sum_{i > 0} \sum_{0<2j \le i}
   {(-n)^k k! \over \lambda^!} \cdot
   {m_i \over (m_{i-j}+1)(m_j+1+\delta_{i, 2j})} \cdot \nonumber \\
& &\quad\quad\quad \cdot
   j(i-j) \cdot {(m_{i-j}+1)(m_j+1+\delta_{i, 2j}) \over 1+ \delta_{i, 2j}}
   \nonumber \\
&=&{(-n)^k k! \over \lambda^!} \sum_{i > 0} m_i \sum_{0 < 2j \le
i}
   {j(i-j) \over 1+ \delta_{i, 2j}}
   = {(-n)^k k! \over \lambda^!} \sum_{i > 0} m_i \left ( {1 \over 2}
   \sum_{0 < j < i} j(i-j) \right )
   \nonumber \\
&=&{(-n)^k k! \over \lambda^!} \sum_{i > 0} m_i \cdot {i^3-i \over
12}
   ={(-n)^k k! \over \lambda^!} \cdot {1 \over 12} \sum_{i > 0} (i^3m_i -im_i).
\end{eqnarray}
Similarly, we see that the total contribution of the second source
with $i < 0$ to the coefficient of $\mathfrak
a_{\lambda}(\tau_{*}(e\alpha))$ in $\mathfrak a_n^{(k+1)}(\alpha)$
is equal to
\begin{eqnarray} \label{derivative2.4}
{(-n)^k k! \over \lambda^!} \cdot {1 \over 12} \sum_{i < 0}
(i^3m_i -im_i).
\end{eqnarray}

Noting $n=\sum_i i m_i$, we see from (\ref{derivative2.2}),
(\ref{derivative2.3}) and (\ref{derivative2.4}) that the
coefficient of the term $\mathfrak a_{\lambda}(\tau_{*}(e\alpha))$
in the derivative $\mathfrak a_n^{(k+1)}(\alpha)$ is:

\begin{eqnarray*}
& &{(-n)^k k! \over \lambda^!} \left ( nc-knc
   + {n \over 12} \sum_i i^2m_i -{1 \over 12} \sum_{i}
  i^3m_i \right) \\
& &+ {(-n)^k k! \over \lambda^!} \cdot {1 \over 12}
   \sum_{i} (i^3m_i - im_i)  \\
&=&{(-n)^k k! \over \lambda^!} (nc-knc-2nc)
   = (-n)^{k+1} (k+1)! {1 - \lambsq \over 24 \lambda^!}.
\end{eqnarray*}

Combining (i) and (ii), we have proved the theorem for $\mathfrak
a_n^{(k+1)}(\alpha)$.
\end{demo}

\begin{remark} \label{derivative3.1}
Let $\alpha\in \mathfrak A_X$. From the proof of
Theorem~\ref{derivative2.1}, we can read that
\begin{eqnarray*}
&& \left ( \sum_{\ell(\lambda) = k+1, |\lambda|=n}
   {1 \over \lambda^!} \mathfrak a_{\lambda}(\tau_{*}\alpha) -
   \sum_{\ell(\lambda) = k-1, |\lambda|=n}
   {\lambsq -1 \over 24\lambda^!}
   \mathfrak a_{\lambda}(\tau_{*}(e\alpha)) \right )^\prime  \\
&=&-n(k+1) \left ( \sum_{\ell(\lambda) = k+2, |\lambda|=n}
   {1 \over \lambda^!} \mathfrak a_{\lambda}(\tau_{*}\alpha) -
   \sum_{\ell(\lambda) = k, |\lambda|=n}
   {\lambsq - 1 \over 24\lambda^!}
   \mathfrak a_{\lambda}(\tau_{*}(e\alpha)) \right ).
\end{eqnarray*}
Note that when $n=0$, this formula is not covered by
Theorem~\ref{derivative2.1}. Applying this formula twice in the
second equality below, we obtain that for a fixed constant $d$,
\begin{eqnarray*}
& &\left ( \sum_{\ell(\lambda) = k+1, |\lambda|=n}
   {1 \over \lambda^!} \mathfrak a_{\lambda}(\tau_{*}\alpha) -
   \sum_{\ell(\lambda) = k-1, |\lambda|=n}
   {\lambsq + d \over 24\lambda^!}
   \mathfrak a_{\lambda}(\tau_{*}(e\alpha)) \right )^\prime  \\
&=&\left ( \sum_{\ell(\lambda) = k+1, |\lambda|=n}
   {1 \over \lambda^!} \mathfrak a_{\lambda}(\tau_{*}\alpha) -
   \sum_{\ell(\lambda) = k-1, |\lambda|=n}
   {\lambsq -1 \over 24\lambda^!}
   \mathfrak a_{\lambda}(\tau_{*}(e\alpha)) \right )^\prime  \\
& &- {d+1 \over 24} \left (
   \sum_{\ell(\lambda) = k-1, |\lambda|=n}
   {1 \over \lambda^!} \mathfrak a_{\lambda}(\tau_{*}(e\alpha))
   \right )^\prime  \\
&=&-n(k+1) \left ( \sum_{\ell(\lambda) = k+2, |\lambda|=n}
   {1 \over \lambda^!} \mathfrak a_{\lambda}(\tau_{*}\alpha) -
   \sum_{\ell(\lambda) = k, |\lambda|=n}
   {\lambsq +d \over 24\lambda^!}
   \mathfrak a_{\lambda}(\tau_{*}(e\alpha)) \right )   \\
& &- {n(d+1) \over 12} \cdot \sum_{\ell(\lambda) = k, |\lambda|=n}
   {1 \over \lambda^!} \mathfrak a_{\lambda}(\tau_{*}(e\alpha)).
\end{eqnarray*}
We shall need this formula in the proof of Theorem \ref{th_chern}.
\end{remark}

Now we describe a formula of $\mathfrak a_n^{(k)}(\alpha)$ for
general $\alpha\in H^*(X)$.

\begin{theorem} \label{deriv_th}
Let $k \ge 0$, $n \in \Z$, and $\alpha\in H^*(X)$. Then,
$\mathfrak a_n^{(k)}(\alpha)$ equals
\begin{eqnarray*}
& &(-n)^k k! \left ( \sum_{\ell(\lambda) = k+1, |\lambda|=n}
   {1 \over \lambda^!} \mathfrak a_{\lambda}(\tau_{*}\alpha) -
   \sum_{\ell(\lambda) = k-1, |\lambda|=n}
   {\lambsq - 1 \over 24 \lambda^!}
   \mathfrak a_{\lambda}(\tau_{*}(e\alpha)) \right )   \\
&+&\sum_{{\epsilon} \in \{\KX, \KX^2\}}  \,\,\,
   \sum_{\ell(\lambda) = k+1-|{\epsilon}|/2, |\lambda|=n}
   {f_\epsilon(\lambda) \over \lambda^!}
   \mathfrak a_{\lambda}(\tau_{*}(\epsilon\alpha))
\end{eqnarray*}
where all the numbers $f_\epsilon(\lambda)$ are independent of $X$
and $\alpha$.
\end{theorem}
\begin{demo}{Proof}
We see from Lemma~\ref{tau_k_tau_{k-1}}~(ii) and (iii) that
\begin{eqnarray} \label{deriv_th1}
\mathfrak a_n^{(k)}(\alpha) =\sum_{{\epsilon} \in \{1_X, e,\KX,
\KX^2\}} \,\,\, \sum_{\ell(\lambda) = k+1-|{\epsilon}|/2,
|\lambda|=n} {f_\epsilon(\lambda) \over \lambda^!} \mathfrak
a_{\lambda}(\tau_{*}(\epsilon\alpha))
\end{eqnarray}
where all the coefficients $f_\epsilon(\lambda)$ are independent
of $X$ and $\alpha$.

To determine $f_{1_X}(\lambda)$ and $f_{e}(\lambda)$, let $X$ be a
$K3$ surface and $\alpha = 1_X$. Then, $\KX \alpha = 0$ but $e
\alpha \ne 0$. By Theorem~\ref{derivative2.1}, $\mathfrak
a_n^{(k)}(\alpha)$ is equal to
\begin{eqnarray}  \label{deriv_th2}
(-n)^k k! \left ( \sum_{\ell(\lambda) = k+1, |\lambda|=n} {1 \over
\lambda^!} \mathfrak a_{\lambda}(\tau_{*}\alpha) -
\sum_{\ell(\lambda) = k-1, |\lambda|=n} {\lambsq -1 \over 24
\lambda^!} \mathfrak a_{\lambda}(\tau_{*}(e\alpha)) \right ).
\end{eqnarray}
Note that the set of all the Heisenberg monomials $\mathfrak
a_{\lambda}(\tau_{*}\alpha), \mathfrak a_{\tilde
\lambda}(\tau_{*}(e\alpha))$ with $\ell(\lambda) = k+1,
|\lambda|=n, \ell(\tilde \lambda) = k-1, |\tilde \lambda|=n$ is
linearly independent. It follows from (\ref{deriv_th1}) and
(\ref{deriv_th2}) that $f_{1_X}(\lambda) = (-n)^k k!$ and
$\displaystyle{f_{e}(\lambda) = -{(-n)^k k! (\lambsq -1) \over
24}}$.
\end{demo}


%
%
%
%
%
%
%
%
\subsection{Formulas for the Chern Character operators}

\begin{lemma} \label{f_{12}}
 \begin{enumerate}
\item[{\rm (i)}] Let $n \ge 1$, $\alpha \in H^*(X)$, and $\mathfrak f \in
\End({\fock})$ with $\mathfrak f'=0$. Then,
\begin{eqnarray*}
[\mathfrak f, \mathfrak a_{-(n+1)}(\alpha)] =-{1 \over n} \cdot \{
[[\mathfrak f, \mathfrak a_{-1}(1_X)]', \mathfrak a_{-n}(\alpha)]
+ [\mathfrak a_{-1}'(1_X), [\mathfrak f, \mathfrak
a_{-n}(\alpha)]]\};
\end{eqnarray*}

\item[{\rm (ii)}] For $i = 1, 2$, let $\mathfrak f_i \in \End({\fock})$ such
that $\mathfrak f_i\cdot |0\rangle=0$ and $\mathfrak f_i'=0$. If
$[\mathfrak f_1, \mathfrak a_{-1}(\alpha)] = [\mathfrak f_2,
\mathfrak a_{-1}(\alpha)]$ for every $\alpha \in H^*(X)$, then
$\mathfrak f_1= \mathfrak f_2$.
 \end{enumerate}
\end{lemma}

\begin{demo}{Proof}
The statement (i) is proved in the Lemma~3.5 in \cite{LQW3}. To
prove (ii), we note that it suffices to show that for every $n \ge
1$ and $\alpha \in H^*(X)$,
\begin{eqnarray} \label{f_{12}1}
[\mathfrak f_1, \mathfrak a_{-n}(\alpha)] = [\mathfrak f_2,
\mathfrak a_{-n}(\alpha)].
\end{eqnarray}
By assumption, (\ref{f_{12}1}) holds for $n=1$. So (\ref{f_{12}1})
follows from (i) and induction.
\end{demo}

Our first main result in this paper is the following theorem.

\begin{theorem} \label{th_chern}
Let $k \ge 0$, and $\alpha\in \mathfrak A_X$. Then, $\mathfrak
G_k(\alpha)$ is equal to
\begin{eqnarray} \label{th_chern1}
- \sum_{\ell(\lambda) = k+2, |\lambda|=0} {1 \over \lambda^!}
\mathfrak a_{\lambda}(\tau_{*}\alpha) + \sum_{\ell(\lambda) = k,
|\lambda|=0} {\lambsq - 2 \over 24\lambda^!} \mathfrak
a_{\lambda}(\tau_{*}(e\alpha)).
\end{eqnarray}
\end{theorem}

\begin{demo}{Proof}
Denote the expression (\ref{th_chern1}) by $\mathfrak f_2$. We
shall apply Lemma~\ref{f_{12}} (ii) to the two operators
$\mathfrak f_1 = \mathfrak G_k(\alpha)$ and $\mathfrak f_2$. Note
that $\mathfrak G_k(\alpha)|0\rangle = 0$ and $\mathfrak
G_k(\alpha)' = 0$. From Theorem~\ref{commutator}~(v) and
Theorem~\ref{derivative2.1}, we see that
\begin{eqnarray*}
& &[\mathfrak G_k(\alpha), \mathfrak a_{-1}(\beta)]
   = {1 \over k!} \cdot \mathfrak a_{-1}^{(k)}(\alpha \beta)   \\
&=&\left ( \sum_{\ell(\lambda) = k+1, |\lambda|=-1}
   {1 \over \lambda^!} \mathfrak a_{\lambda}(\tau_{*}(\alpha \beta)) -
   \sum_{\ell(\lambda) = k-1, |\lambda|=-1}
   {\lambsq -1 \over 24 \lambda^!}
   \mathfrak a_{\lambda}(\tau_{*}(e\alpha\beta)) \right ).
\end{eqnarray*}

Next we need to check the corresponding properties for $\mathfrak
f_2$. Indeed, $\mathfrak f_2 \cdot |0\rangle = 0$ from the
definition of $\mathfrak f_2$. Also, $\mathfrak f_2' = 0$ by
Remark~\ref{derivative3.1}. By Lemma~\ref{tau_k_tau_{k-1}} (i), we
obtain
\begin{eqnarray*}
& &[\mathfrak f_2, \mathfrak a_{-1}(\beta)]  \\
&=&\sum_{\ell(\lambda) = k+1, |\lambda|=-1}
   {1 \over \lambda^!} \mathfrak a_{\lambda}(\tau_{*}(\alpha\beta))
   - \sum_{\ell(\lambda) = k, |\lambda|=0}
   {\lambsq - 2 \over 24(\lambda_1)^!}
   \mathfrak a_{\lambda_1}(\tau_{*}(e\alpha\beta))   \\
&=&\sum_{\ell(\lambda) = k+1, |\lambda|=-1}
   {1 \over \lambda^!} \mathfrak a_{\lambda}(\tau_{*}(\alpha\beta))
   - \sum_{\ell(\lambda_1) = k-1, |\lambda_1|=-1}
   {s(\lambda_1) -1 \over 24(\lambda_1)^!}
   \mathfrak a_{\lambda_1}(\tau_{*}(e\alpha\beta))   \\
&=&[\mathfrak G_k(\alpha), \mathfrak a_{-1}(\beta)]
\end{eqnarray*}
where for a generalized partition $\lambda$ with $\ell(\lambda) =
k$, the generalized partition $\lambda_1$ is obtained from
$\lambda$ by subtracting $1$ from the multiplicity of $1$.

Thus, we conclude from Lemma~\ref{f_{12}} (ii) that $\mathfrak
G_k(\alpha) = \mathfrak f_2$.
\end{demo}

Now we present a formula of $\mathfrak G_k(\alpha)$ for a general
$\alpha \in H^*(X)$.

\begin{theorem} \label{char_th}
Let $k \ge 0$ and $\alpha\in H^*(X)$. Then, $\mathfrak
G_k(\alpha)$ is equal to
\begin{eqnarray*}
& &- \sum_{\ell(\lambda) = k+2, |\lambda|=0}
   {1 \over \lambda^!} \mathfrak a_{\lambda}(\tau_{*}\alpha)
   + \sum_{\ell(\lambda) = k, |\lambda|=0}
   {\lambsq - 2 \over 24\lambda^!}
   \mathfrak a_{\lambda}(\tau_{*}(e\alpha))  \\
&+&\sum_{{\epsilon} \in \{\KX, \KX^2\}}  \,\,\,
   \sum_{\ell(\lambda) = k+2-|{\epsilon}|/2, |\lambda|=0}
   {g_\epsilon(\lambda) \over \lambda^!}
   \mathfrak a_{\lambda}(\tau_{*}(\epsilon\alpha))
\end{eqnarray*}
where all the numbers $g_\epsilon(\lambda)$ are independent of $X$
and $\alpha$.
\end{theorem}

\begin{demo}{Proof}
First of all, by the formula (5.6) in \cite{LQW1}, we have
\begin{eqnarray*}
\mathfrak G_0(\alpha) = - \mathfrak L_0(\alpha) = -
\sum_{\ell(\lambda) = 2, |\lambda|=0} {1 \over \lambda^!}
\mathfrak a_{\lambda}(\tau_{*}\alpha).
\end{eqnarray*}
So our theorem is true for $k = 0$. In the following, we assume $k
> 0$ and $\alpha \ne 0$.

Next, we write the operator $\mathfrak G_k(\alpha) \in
\End(\fock)$ as a linear combination of Heisenberg monomials
$\displaystyle{\prod_i (\mathfrak a_i(\alpha_{i,1}) \cdots
\mathfrak a_i(\alpha_{i, m_i}))}$. For degree reasons, $\sum_i i
m_i = 0$ and
\begin{eqnarray}   \label{char_th1}
-2 \sum_i m_i + \sum_{i,j} |\alpha_{i, j}|= 2k+|\alpha|.
\end{eqnarray}
In other words, we regard $\mathfrak G_k(\alpha)$ as an element in
the completion of the universal enveloping algebra of the
Heisenberg algebra. We can always do so because $\fock$ is
irreducible as a representation of the Heisenberg algebra. Now
rewrite $\mathfrak G_k(\alpha)$ as
\begin{eqnarray}  \label{char_th2}
\mathfrak G_k(\alpha) = \sum_{|\lambda| = 0} \mathfrak g_\lambda
\end{eqnarray}
where for each fixed generalized partition $\lambda = (\cdots
(-2)^{m_{-2}}(-1)^{m_{-1}} 1^{m_1}2^{m_2} \cdots)$ with $|\lambda|
= 0$, the operator $\mathfrak g_\lambda$ stands for the component
in $\mathfrak G_k(\alpha)$ containing all the expressions of the
form $\displaystyle{\prod_i (\mathfrak a_i(\alpha_{i,1}) \cdots
\mathfrak a_i(\alpha_{i, m_i}))}$.

Fix $n \in \Z$ and $\beta \in H^*(X)$. By
Lemma~\ref{tau_k_tau_{k-1}} (iii) and the Lemma~5.1 in
\cite{LQW3},
\begin{eqnarray}   \label{char_th3}
[\mathfrak G_k(\alpha), \mathfrak a_n(\beta)] = \sum_{{\epsilon}
\in \{1_X, \KX, \KX^2, e\}} \, \, \, \sum_{\ell(\mu) =
k+1-|\epsilon|/2, |\mu| = n} d(\epsilon, \mu) \mathfrak
a_{\mu}(\tau_*(\epsilon\alpha\beta))
\end{eqnarray}
where all the coefficients $d(\epsilon, \mu) \in \mathbb Q$ are
independent of the surface $X$ and the cohomology classes $\alpha,
\beta \in H^*(X)$. It follows that the generalized partitions
$\lambda$ in (\ref{char_th2}) must satisfy $k \le \ell(\lambda)
\le (k+2)$. Moreover, for a fixed generalized partition $\lambda$
with $|\lambda| = 0$ and $k \le \ell(\lambda) \le (k+2)$, and for
every $i$ with $m_{-i} > 0$, we have
\begin{eqnarray}  \label{char_th4}
&[\mathfrak g_\lambda, \mathfrak a_i(\beta)] =
   \left\{
      \everymath{\displaystyle}
      \begin{array}{ll}
         d(1_X, \lambda_{-i}) \mathfrak a_{\lambda_{-i}}(\tau_*(\alpha\beta)),
           & \ell(\lambda) = k+2     \\
         d(\KX, \lambda_{-i}) \mathfrak a_{\lambda_{-i}}(\tau_*(\KX\alpha\beta)),
           & \ell(\lambda) = k+1     \\
         \mathfrak a_{\lambda_{-i}}(\tau_*((d(e, \lambda_{-i})e
             +d(\KX^2, \lambda_{-i})\KX^2)\alpha\beta)),
           & \ell(\lambda) = k.
      \end{array}
    \right.
\end{eqnarray}
where $\lambda_{-i}$ is the generalized partition obtained from
$\lambda$ by subtracting $1$ from the multiplicity of $(-i)$. Note
that $d(\epsilon, \lambda_{-i})$ is independent of $X$, $\alpha$
and $\beta$.

By choosing special surfaces $X$ and $\alpha, \beta \in H^*(X)$,
we now claim that
$${d(1_X, \lambda_{-i}) \over i m_{-i}}, \,\, {d(\KX, \lambda_{-i}) \over
i m_{-i}}, \,\, {d(e, \lambda_{-i}) \over i m_{-i}} e + 
{d(\KX^2, \lambda_{-i}) \over i m_{-i}}\KX^2  \eqno (4.13)$$
are independent of $i$ as long as $m_{-i} > 0$. Indeed, to prove that
the 3rd item in (4.13) is independent of $i$ as long as
$m_{-i} > 0$, we may assume that $\mathfrak g_\lambda \ne 0$ 
by (4.12). Note that $|\lambda| = 0$ and
$\ell(\lambda) = k$. By (4.9), $\sum_{i,j} |\alpha_{i,
j}|= 4k+|\alpha|$ for every Heisenberg monomials
$\displaystyle{\prod_i (\mathfrak a_i(\alpha_{i,1}) \cdots
\mathfrak a_i(\alpha_{i, m_i}))}$ contained in $\mathfrak
g_\lambda$. On the other hand, $|\alpha_{i, j}| \le 4$
and thus $\sum_{i,j} |\alpha_{i, j}| \le 4\ell(\lambda) = 4k$. 
It follows that $|\alpha| = 0$ and $|\alpha_{i, j}| = 4$ for
every $i$ and $j$. So
$\mathfrak g_\lambda = b \cdot \mathfrak a_{\lambda}(\tau_*[x])
$
where $b$ is a nonzero number, and $[x] \in H^4(X)$ is the
cohomology class corresponding to a point $x \in X$. 
Since $\alpha \ne 0$ and $|\alpha| = 0$, $\alpha = \tilde b
\cdot 1_X$ for some nonzero number $\tilde b$. 
Choosing $\beta = 1_X$, we see from (4.12) that
\begin{eqnarray*}
m_{-i} i b \cdot \mathfrak a_{\lambda_{-i}}(\tau_*[x]) = [\mathfrak
g_\lambda, \mathfrak a_i(1_X)] = \tilde b \cdot \mathfrak
a_{\lambda_{-i}}(\tau_*(d(e, \lambda_{-i})e
             +d(\KX^2, \lambda_{-i})\KX^2)).
\end{eqnarray*}
So $b[x] = \tilde b \cdot (d(e, \lambda_{-i})e +d(\KX^2,
\lambda_{-i})\KX^2)/(i m_{-i})$. Thus
$\displaystyle{{d(e, \lambda_{-i}) \over i m_{-i}} e + {d(\KX^2,
\lambda_{-i}) \over i m_{-i}}\KX^2}$
is independent of $i$ as long as $m_{-i} > 0$. Similarly, by
choosing $\alpha \in H^2(X)$ with $\KX \alpha \ne 0$ (resp.
$\alpha \in H^4(X)-\{0\}$), we see that ${d(\KX, \lambda_{-i})/(i m_{-i})}$
(resp. ${d(1_X, \lambda_{-i})/(i m_{-i})}$) is independent of $i$ as long
as $m_{-i} > 0$. This proves (4.13).

Now, for an arbitrary $X$ and $\alpha \in H^*(X)$, we claim that
$\mathfrak G_k(\alpha)$ is equal to
\begin{eqnarray*}
   \mathfrak f
&\stackrel{\rm def}{=}&
   \sum_{{\epsilon} \in \{1_X, \KX, \KX^2, e\}} \, \, \,
   \sum_{\ell(\lambda) = k+2-|\epsilon|/2, |\lambda| = 0}
   {d(\epsilon, \lambda_{-i}) \over i m_{-i}}
   \mathfrak a_{\lambda}(\tau_*(\epsilon\alpha)).
\end{eqnarray*}
%
Indeed, for every $n \in \Z-\{0\}$ and $\beta \in H^*(X)$, we have
\begin{eqnarray*}
   [\mathfrak f, \mathfrak a_n(\beta)]
&=&\sum_{{\epsilon} \in \{1_X, \KX, \KX^2, e\}} \, \, \,
   \sum_{\ell(\lambda) = k+2-|\epsilon|/2, |\lambda| = 0}
   n m_{-n} {d(\epsilon, \lambda_{-i}) \over i m_{-i}}  \cdot
   \mathfrak a_{\lambda_{-n}}(\tau_*(\epsilon\alpha\beta))  \\
&=&\sum_{{\epsilon} \in \{1_X, \KX, \KX^2, e\}} \, \, \,
   \sum_{\ell(\lambda) = k+2-|\epsilon|/2, |\lambda| = 0}
   d(\epsilon, \lambda_{-n})  \cdot
   \mathfrak a_{\lambda_{-n}}(\tau_*(\epsilon\alpha\beta))  \\
&=&[\mathfrak G_k(\alpha), \mathfrak a_n(\beta)]
\end{eqnarray*}
where we have used (4.13) in the second equality, and
(4.10) and (4.12) in the last equality. 
Since $\fock$ is irreducible, we see from Schur's lemma that
$(\mathfrak G_k(\alpha)-\mathfrak f)$ must be a scalar multiple of
the identity operator. Recall that the bidegree of $(\mathfrak
G_k(\alpha)-\mathfrak f)$ is $(0, 2k+|\alpha|)$ which is
nontrivial since $k > 0$. Therefore, $(\mathfrak
G_k(\alpha)-\mathfrak f) = 0$. So
\begin{eqnarray*}
\mathfrak G_k(\alpha)
 = \mathfrak f
=  \sum_{{\epsilon} \in \{1_X, \KX, \KX^2, e\}} \, \, \,
   \sum_{\ell(\lambda) = k+2-|\epsilon|/2, |\lambda| = 0}
   {d(\epsilon, \lambda_{-i}) \over i m_{-i}}
   \mathfrak a_{\lambda}(\tau_*(\epsilon\alpha)).
\end{eqnarray*}

Finally, since the numbers $d(\epsilon, \lambda_{-i})/(i m_{-i})$ 
are independent of $X$ and $\alpha \in H^*(X)$, 
by using Theorem~4.6 and 
an argument similar to that in the proof of Theorem~4.4, 
we conclude that $d(1_X, \lambda_{-i})/(i m_{-i}) = -1/\lambda^!$ and
$d(e, \lambda_{-i})/(i m_{-i}) = (\lambsq - 2)/(24\lambda^!)$. In particular,
$d(e, \lambda_{-i})/(i m_{-i})$ is independent of $i$ whenever $m_{-i} > 0$.
So by (4.13), $d(K^2, \lambda_{-i})/(i m_{-i})$ is also 
independent of $i$ whenever $m_{-i} > 0$. Now denoting 
$d(\epsilon, \lambda_{-i})/(i m_{-i})$ by $g_{\epsilon}(\lambda)/\lambda^!$
for $\epsilon = K$ and $K^2$, we have proved the theorem.

\end{demo}

%

\begin{corollary} \label{g_{1_X, e}}
Let $n \ge 1$, $k \ge 0$, and $\alpha \in H^*(X)$. Then, $\quad
G_k(\alpha, n)$ is equal to
\begin{eqnarray*}
& &\sum_{0 \le j \le k} \sum_{\lambda
   \vdash (j+1) \atop \ell(\lambda)=k-j+1}
   {(-1)^{|\lambda|-1} \over \lambda^! \cdot |\lambda|!}
   \cdot {\bf 1}_{-(n-j-1)} \mathfrak a_{-\lambda}(\tau_*\alpha)|0\rangle \\
&+&\sum_{0 \le j \le k} \sum_{\lambda \vdash (j+1) \atop
\ell(\lambda)=k-j-1}
   {(-1)^{|\lambda|} \over \lambda^! \cdot |\lambda|!}
   \cdot {|\lambda| + \lambsq - 2 \over 24}
   \cdot {\bf 1}_{-(n-j-1)}
   \mathfrak a_{-\lambda}(\tau_*(e\alpha))|0\rangle  \\
&+&\sum_{\epsilon \in \{\KX, \KX^2\} \atop 0 \le j \le k}
   \sum_{\lambda \vdash (j+1) \atop
\ell(\lambda)=k-j+1-|{\epsilon}|/2}
   {(-1)^{|\lambda|}g_{\epsilon}(\lambda+(1^{j+1}))
   \over \lambda^! \cdot |\lambda|!} \cdot {\bf 1}_{-(n-j-1)}
   \mathfrak a_{-\lambda}(\tau_*(\epsilon\alpha))|0\rangle
\end{eqnarray*}
where ${\bf 1}_{-(n-j-1)}$ denotes $\mathfrak
a_{-1}(1_X)^{n-j-1}/(n-j-1)!$ when $(n-j-1) \ge 0$ and is $0$ when
$(n-j-1) < 0$, the universal function $g_{\epsilon}$ is from
Theorem~\ref{char_th}, and $\lambda+(1^{j+1})$ is the partition
obtained from $\lambda$ by adding $(j+1)$ to the multiplicity of
$1$.
\end{corollary}

\begin{demo}{Proof}
Let $1_{\Xn}$ be the fundamental class of the Hilbert scheme
$\Xn$. Then, we have $1_{\Xn} = \displaystyle{{1 \over n!} \cdot
\mathfrak a_{-1}(1_X)^n}|0\rangle$. From the definition of Chern
character operators, we see that $\mathfrak G_k(\alpha)|0\rangle =
0$ and $G_k(\alpha, n) = \mathfrak G_k(\alpha) 1_{\Xn}$. Also, by
Theorem~\ref{char_th} and Lemma~\ref{tau_k_tau_{k-1}}~(i), we have
$[\ldots [\mathfrak G_k(\alpha), \mathfrak a_{n_1}(\alpha_1)],
\ldots], \mathfrak a_{n_{k+2}}(\alpha_{k+2})]=0$ whenever $n_1,
\ldots, n_{k+2} < 0$. Combining all of these observations, we
obtain
\begin{eqnarray*}
 G_k(\alpha, n)
 &=&{1 \over n!} \cdot \mathfrak G_k(\alpha)
    \mathfrak a_{-1}(1_X)^n|0\rangle  \\
 &=& {1 \over n!} \cdot \sum_{u=1}^{k+1} {n \choose u}
    \mathfrak a_{-1}(1_X)^{n-u}
    [\cdots [\mathfrak G_k(\alpha), \underbrace{\mathfrak a_{-1}(1_X)],
    \cdots], \mathfrak a_{-1}(1_X)]}_{u ~\rm{times}}|0\rangle  \\
 &=&\sum_{j=0}^{k} {1 \over (j+1)!} \cdot {\bf 1}_{-(n-j-1)}
    [\cdots [\mathfrak G_k(\alpha), \underbrace{\mathfrak a_{-1}(1_X)],
    \cdots], \mathfrak a_{-1}(1_X)]}_{j+1 ~\rm{times}}|0\rangle .
\end{eqnarray*}
Now our corollary follows immediately from Theorem~\ref{char_th}.
\end{demo}

\begin{remark}
(i) Let $X$ be a projective surface with numerically trivial
canonical class. Recall from Sect.~\ref{sect_general} that the
classes $G_k(\alpha,n)$, where $0\le k<n$ and $\alpha$ runs over a
linear basis of $H^*(X)$, form a set of ring generators for the
cohomology ring $H^*(\Xn)$. Rewrite the cup product of $k$ ring
generators as
\begin{eqnarray} \label{eq_operator}
&&G_{i_1}(\alpha_1,n)\cdot G_{i_2}(\alpha_2,n)
 \cdot ... \cdot G_{i_k}(\alpha_k,n)  \nonumber  \\
&=&\mathfrak G_{i_1}(\alpha_1)\mathfrak G_{i_2}(\alpha_1)
  \cdots \mathfrak G_{i_k}(\alpha_k)
  \cdot {1 \over n!} \mathfrak a_{-1}(1_X)^n \vac.
\end{eqnarray}
Applying Theorem~\ref{th_chern} to $\mathfrak G_{i_j}(\alpha_j)$
and using Heisenberg algebra commutation relation, we can express
(\ref{eq_operator}) as a linear combination of Heisenberg
monomials. Such a procedure is analogous to the one in the proof
of Corollary~\ref{g_{1_X, e}}. In this way, Theorem~\ref{th_chern}
provides us a complete description of the ring $H^*(\Xn)$. There
has also been another very different description of the ring
$H^*(\Xn)$ in \cite{LS2}.

(ii) In addition, when specialized to surfaces with numerically
trivial canonical classes, our Corollary~\ref{g_{1_X, e}}
establishes a (slightly corrected) conjectural formula in
Remark~4.15, \cite{LS2}. (A factor $2$ was missing in the term
involving the Euler class $e$ in the conjectural formula of
Lehn-Sorger).

(iii) For a general projective surface $X$, Theorem~\ref{char_th}
provides a partial description to the cohomology ring of the
Hilbert scheme $\Xn$.
\end{remark}

\begin{remark} \label{rmk-g_{1_X, e}}
{\rm Let $[x] \in H^4(X)$ be the cohomology class corresponding to
a point $x \in X$. Assume $\sum\limits_{i=1}^s (k_i+2) = 2n$. By
Theorem \ref{commutator} (v), $[\mathfrak G_{k_i}([x]), \mathfrak
a_{-1}([x])] = 0$ for all $1 \le i \le s$. Using
Theorem~\ref{char_th} and Corollary~\ref{g_{1_X, e}} repeatedly,
we see that the intersection number $\displaystyle{\prod_{i=1}^s
G_{k_i}([x], n)} \in H^{4n}(\Xn) \cong \mathbb Q$ is equal to
\begin{eqnarray*}
\sum_{0 \le j_1 \le k_1, \ldots,  0 \le j_s \le k_s \atop (j_1+1)+
\ldots + (j_s+1) = n} \prod_{i=1}^s \sum_{ \lambda_i \vdash
(j_i+1) \atop \ell(\lambda_i)=k_i-j_i+1} {(-1)^{|\lambda_i|-1}
\over (\lambda_i)^! |\lambda_i|!}.
\end{eqnarray*}
Note that this number is independent of the projective surface $X$
(see \cite{LQW3}).}
\end{remark}
\section{\bf Action of the $\mathcal W$ algebras on $\fock$}
\label{sect_walghilb}

In this section, we formulate the main theorems on connections
between $\mathcal W$ algebras and Hilbert schemes. The proof of
Theorem \ref{th_main}, which is very technical and uses the
operator product expansion technique, will be postponed to
Sect.~\ref{sect_proof}.

First of all, we introduce the following definitions.
\begin{definition} \label{J}
Let $X$ be a smooth projective surface.
\begin{enumerate}
\item[{\rm (i)}] For $p \ge 0$, $n \in \Z$ and $\alpha \in H^*(X)$,
define $\jj^p_n(\alpha) \in \End(\fock)$ to be
\begin{eqnarray*}
  p! \cdot\left( -\sum_{\ell(\lambda) = p+1, |\lambda|=n}
  \frac{1}{\lambda^!} \mathfrak a_{\lambda}(\tau_{*}\alpha)
    + \sum_{\ell(\lambda) = p-1, |\lambda|= n} \frac{\lambsq
     + n^2 - 2}{24\lambda^!} \mathfrak
     a_{\lambda}(\tau_{*}(e\alpha))
\right);
\end{eqnarray*}

\item[{\rm (ii)}] We define $\Wax$ to be the linear span of the
identity operator ${\rm Id}_{\fock}$ and the operators
$\jj^p_n(\alpha)$ in $\End(\fock)$, where $p \ge 0, n \in \Z$ and
$\alpha \in H^*(X)$.
\end{enumerate}
\end{definition}

Some of the operators $\jj^p_n(\alpha)$ can be identified with the
familiar ones, by using the above definition,
Theorem~\ref{th_chern} and Theorem~\ref{derivative2.1}. For
example, for $\alpha \in H^*(X)$, we see that $\jj^0_n(\alpha) =
-\mathfrak a_{n}(\alpha)$ and $\jj^1_n(\alpha) = \mathfrak
L_{n}(\alpha)$. For $\alpha \in \mathfrak A_X$, we obtain
$\jj^{p}_0(\alpha) = p! \cdot\mathfrak G_{p-1}(\alpha)$ and
$\jj^{p}_{-1}(\alpha)= - \mathfrak a_{-1}^{(p)}(\alpha)$. In
general, we have the following.

\begin{lemma} \label{lem_geom}
Given $p \ge 0$, $\alpha\in \mathfrak A_X$ and $\beta \in H^*(X)$,
we have
\begin{eqnarray*}
[\mathfrak G_p(\alpha), \mathfrak a_{n}(\beta)] =\frac{n}{p!}
\cdot \jj^p_n(\alpha\beta).
\end{eqnarray*}
\end{lemma}
\begin{demo}{Proof}
Follows from Definition~\ref{J}~(i), Theorem~\ref{th_chern}, and
Lemma \ref{tau_k_tau_{k-1}} (i).
\end{demo}

By Lemma~\ref{lem_geom}, we may regard the linear operators
$\jj^p_n(\alpha) \in \End(\fock)$ being geometric since both the
Chern character operators $\mathfrak G_p(\alpha)$ and the
Heisenberg operators $\mathfrak a_{n}(\beta)$ are geometric by the
constructions in Sect.~\ref{sect_general}.

We are interested in the commutation relation among the operators
$\jj^p_m (\alpha)$. To this end, it is useful to adopt the vertex
algebra language, cf.  \cite{Bor, FB, Kac}, to our present setup.
Our convention of {\it fields} is to write them in a form
\begin{eqnarray*}
\phi(z) =\sum_n \phi_n z^{-n-\Delta}
\end{eqnarray*}
where $\Delta$ is the conformal weight of the field $\phi (z)$.
Define
\begin{eqnarray*}
\partial \phi(z) =\sum_n (-n-\Delta) \phi_n z^{-n-\Delta-1}
\end{eqnarray*}
which is called the {\it derivative field} of $\phi(z)$. We set
$\phi_-(z)= \sum_{n \ge 0} \phi_{n} z^{-n-\Delta}$ and $\phi_+(z)
= \sum_{n < 0} \phi_{n} z^{-n-\Delta}$.
If $\psi(z)$ is another vertex operator, we define a new vertex
operator, which is called the {\it normally ordered product} of
$\phi(z)$ and $\psi(z)$, to be:
\begin{eqnarray*}
:\phi(z)\psi(z): \,\,\, = \,\,\, \phi_+(z) \psi(z) +
(-1)^{\phi\psi} \psi(z) \phi_-(z)
\end{eqnarray*}
where $(-1)^{\phi\psi}$ is $-1$ if both $\phi(z)$ and $\psi(z)$
are odd (i.e. fermionic) fields and $1$ otherwise. Inductively we
can define the normally ordered product of $k$ vertex operators
$\phi_1(z), \phi_2(z), \ldots, \phi_k(z)$ from right to left by
\begin{eqnarray*}
:\phi_1(z) \phi_2(z) \cdots \phi_k(z): \,\,\, = \,\,\, : \phi_1(z)
(:\phi_2(z) \cdots \phi_k(z):):.
\end{eqnarray*}

For $\alpha \in H^*(X)$, we define a vertex operator (i.e. a
field) $a(\alpha)(z)$ by putting
\begin{eqnarray*}
a(\alpha)(z) = \sum_{n \in \Z} \mathfrak a_{n}(\alpha) z^{-n-1}.
\end{eqnarray*}
The field $:a(z)^p:(\tau_{*}\alpha)$ is defined to be $\sum_i
:a(\alpha_{i,1})(z) a(\alpha_{i,2})(z) \cdots a(\alpha_{i,p})(z):$
if we write $\tau_{p*}(\alpha) =\sum_i \alpha_{i,1}\otimes
\alpha_{i,1}\otimes \ldots \otimes \alpha_{i,p} \in
H^*(X)^{\otimes p}$. Note that $:a(z)^{p}: (\ta)$ is of conformal
weight $p$. So we rewrite $:a(z)^{p}:(\ta)$ componentwise as
\begin{eqnarray*}
:a(z)^{p}: (\ta)= \sum_m :a^{p}:_m (\ta)\; z^{-m-p},
\end{eqnarray*}
where $:a^{p}:_m (\ta) \in \End(\fock)$ is the coefficient of
$z^{-m-p}$ (i.e. the $m$-th Fourier component of the field
$:a(z)^{p}: (\ta)$), and maps $H^*(X^{[n]})$ to $H^*(X^{[n+m]})$.
Similarly, for $r \ge 1$, we can define the field $:(\partial^r
a(z)) a(z)^{p-1}:(\tau_{*}\alpha)$, and define the operator
$:(\partial^r a)a^{p-1}:_m(\tau_{*}\alpha)$ as the coefficient of
$z^{-m- r-p}$ in $:(\partial^r a(z))a(z)^{p-1}:(\tau_{*}\alpha)$.

We remark that when the variable $z$ or the cohomology class
$\alpha$ is clear or irrelevant in the context, we shall drop $z$
or $\alpha$ from the notations of fields. For instance, we shall
sometimes use $:a^{p}: (\ta)$ or $:a^{p}:$, $:(\partial^r
a)a^p:(\tau_{*}\alpha)$ or $:(\partial^r a)a^p:$, etc to stand for
$:a(z)^{p}: (\ta)$, $:(\partial^r a(z))a(z)^p:(\tau_{*}\alpha)$,
etc respectively.

\begin{lemma} \label{lem_fourier}
In terms of fields, the operator $\jj^p_m(\alpha)$ can be
rewritten as:
\begin{eqnarray} \label{eq_field}
 & & -\frac1{(p+1)} :a^{p+1}:_m (\ta) \nonumber + \frac1{24}
                 p(m^2-3m-2p) :a^{p-1}:_m(\tea) \nonumber \\
 & &+ \frac1{24} p(p-1) :(\pa^2 a)\, a^{p-2}:_m (\tea).
 \end{eqnarray}
\end{lemma}

\begin{demo}{Proof}
First of all, by the definition of $:a^{p+1}:_m (\ta)$, we have
\begin{eqnarray*}
 p! \cdot \sum_{\ell(\lambda) = p+1, |\lambda|=m}
  \frac{1}{\lambda^!} \mathfrak a_{\lambda}(\tau_{*}\alpha)
 &=& \frac1{(p+1)} \sum_{i_1 + \ldots +i_{p+1}=m} :\mathfrak a_{i_1}
     \cdots \mathfrak a_{i_{p+1}} :(\ta) \\
 &=&\frac1{(p+1)} :a^{p+1}:_m (\ta).
\end{eqnarray*}

Next, recalling the definition of $\lambsq$, we obtain
\begin{eqnarray*}
 && p! \cdot \sum_{\ell(\lambda) = p-1, |\lambda|= m}
 \frac{\lambsq + m^2 - 2}{\lambda^!} \mathfrak
     a_{\lambda}(\tau_{*}(e\alpha))  \\
 &=& p\cdot \sum_{i_1 + \ldots +i_{p-1}=m}
     \left(\sum_{b=1}^{p-1}i_b^2 +m^2-2
     \right)
    :\mathfrak a_{i_1}\cdots \mathfrak a_{i_{p-1}}:(\tau_{*}(e\alpha)) \\
 &=& p\cdot \sum_{i_1 + \ldots +i_{p-1}=m}
     \left(\sum_{b=1}^{p-1}(i_b+1)(i_b+2)+m^2-3m-2p
     \right)
     :\mathfrak a_{i_1} \cdots \mathfrak a_{i_{p-1}}:(\tau_{*}(e\alpha)) \\
 &=& p(p-1) :(\pa^2 a)\, a^{p-2}:_m (\tea)
     + p(m^2-3m-2p) :a^{p-1} :_m(\tea).
\end{eqnarray*}

Now the lemma follows from the definition of the operator
$\jj^p_m(\alpha)$.
\end{demo}

\begin{remark} \rm
More explicitly, we can rewrite the operator $\jj^p_m(\alpha)$ as
the $m$-th Fourier component of a vertex operator as follows:
\begin{eqnarray*}
\jj^p_m(\alpha)
 &=& -\frac1{(p+1)} :a^{p+1} :_m (\ta)
 + \frac{p}{24}  (\pa^2:a^{p-1}:)_m(\tea)   \\
 && + \frac{(p+1)p}{12} (\pa :a^{p-1} :)_m(\tea)
 + \frac{p(p^2-p-2)}{24}  :a^{p-1} :_m(\tea)   \\
 & &+\frac1{24} p(p-1) :(\pa^2 a) a^{p-2}:_m (\tea).
\end{eqnarray*}
In practice, (\ref{eq_field}) suffices for the purpose of
computations in the next section. Also, since $\jj^{p+1}_0(\alpha)
= (p+1)! \cdot \mathfrak G_{p}(\alpha)$ for $\alpha \in \mathfrak
A_X$, we conclude that the Chern character operator $\mathfrak
G_{p}(\alpha)$ with $\alpha \in \mathfrak A_X$ is the zero-mode of
a vertex operator.
\end{remark}


One of the main results in this paper is about the commutator
$[\jj^p_m(\alpha), \jj^q_n(\beta)]$. To state the formula, we
define an integer $\Omega_{m,n}^{p,q}$ for $m,n,p,q \in \Z$ as
follows:
\begin{eqnarray}
 \Omega_{m,n}^{p,q}
  &=& mp^3n^2 +3mp^2n^2q -p^2nq +p^2qn^3-3mp^2n^2 +pnq  \nonumber \\
  && +3m^2pnq -3mpn^2q -m^3q^2p -pqn^3 -mpq +m^3pq  \nonumber \\
  && +mpq^2 +2mpn^2 -3m^2pnq^2 -2m^2nq +3m^2nq^2 -m^2nq^3 .
\label{eq_number}
\end{eqnarray}

\begin{theorem}  \label{th_main}
The vector space $\Wax$ is closed under the Lie bracket. More
explicitly, for $m,n \in \Z$, and $\alpha, \beta \in H^*(X)$, we
have
\begin{eqnarray*}
 [\jj^p_m(\alpha), \jj^q_n(\beta)]
 &=& (q m -p n) \cdot \jj^{p+q-1}_{m+n} (\alpha\beta)
  - \frac{\Omega_{m,n}^{p,q}}{12} \cdot \jj^{p+q-3}_{m+n}(e\alpha\beta)
\end{eqnarray*}
where $(p,q) \in \Z_+^2$ except for the {\em unordered} pairs
$(0,0), (1,0), (2,0)$ and $(1,1)$. In addition, for these four
exceptional cases, we have
\begin{eqnarray*}
 {[} \jj^0_m(\alpha), \jj^0_n(\beta)]
 &=& - m \delta_{m,-n} \int_X(\alpha\beta) \cdot {\rm Id}_{\fock}, \\
 {[} \jj^1_m(\alpha), \jj^0_n(\beta)]
 &=& - n \cdot \jj^0_{m+n} (\alpha\beta),\\
 {[} \jj^2_m(\alpha), \jj^0_n(\beta)]
 &=& - 2n \cdot \jj^1_{m+n} (\alpha\beta)
     + \frac{m^3-m}{6}\delta_{m,-n} \int_X(e\alpha\beta)
     \cdot {\rm Id}_{\fock},\\
 {[} \jj^1_m(\alpha), \jj^1_n(\beta)]
 &=& (m-n) \cdot \jj^1_{m+n} (\alpha\beta)
     + \frac{m^3-m}{12}\delta_{m,-n} \int_X(e\alpha\beta)
     \cdot {\rm Id}_{\fock}.
\end{eqnarray*}
\end{theorem}
The proof of this theorem, which uses the operator product
expansion (OPE) method in the theory of vertex algebras with some
appropriate modifications, will be given in Sect.
\ref{sect_proof}.
Note that our algebra $\Wax$ contains as subalgebras the
Heisenberg algebra of Nakajima-Grojnowski generated by the
operators $\jj^0_m(\alpha) =-\mathfrak a_m(\alpha)$ and the
Virasoro algebra of Lehn generated by the operators
$\jj^1_m(\alpha)=\mathfrak L_m(\alpha)$.

\begin{remark}
Since $\mathfrak d = \mathfrak G_1(1_X) \equiv \frac12
\jj^2_0(1_X)$ modulo the $\KX$-term (cf. Theorem~\ref{char_th},
also cf. \cite{LQW1, FW}), a special case of Theorem~\ref{th_main}
reads that
\begin{eqnarray*}
(\jj^p_n(\alpha))' = \frac12 [\jj^2_0(1_X), \jj^p_n(\alpha)] = -n
\jj^{p+1}_n(\alpha) - \frac{(n^3-n)p}{12} \jj^{p-1}_n(e\alpha)
\end{eqnarray*}
for $p\ge 1$ and $\alpha\in \mathfrak A_X$. Note that setting
$d=n^2-2$ in Remark~\ref{derivative3.1} provides a totally
different way of proving the same formula.
\end{remark}

We denote by $\mathfrak B_X =\{\alpha\in H^*(X)|e \alpha = \KX
\alpha = 0\}$ which is an ideal in the ring $H^*(X)$. Obviously,
$\mathfrak B_X \subset \mathfrak A_X$. Denote by $\Wbx$ the linear
span of $\jj^p_n(\alpha)$, where $p \ge 0, n \in \Z$ and $\alpha
\in \mathfrak B_X$. From the commutation relation in
Theorem~\ref{th_main} we observe that $\Wbx$ is a Lie subalgebra
of the Lie (super)algebra $\Wax$. Recall the (super)algebra
$\widehat{\W}(A)$ introduced in Section~\ref{sect_walg} for a
general ring $A$. The next theorem follows from comparing the
commutator (\ref{eq_gradedcomm}) and the ones given in
Theorem~\ref{th_main}.

\begin{theorem}
Let $X$ be a smooth projective surface. Then, the Lie
(super)algebra $\widehat{\W}(\mathfrak B_X)$ is isomorphic to the
Lie (super)algebra $\Wbx$ by sending $C \mapsto {\rm Id}_{\fock}$
and $\mathfrak L_n^p(\alpha)\mapsto \jj_n^p(\alpha)$, where $p \ge
0$, $n\in \Z$, and $\alpha \in \mathfrak B_X$.
\end{theorem}

\begin{remark} \rm
(i) For a smooth projective surface $X$ with numerically trivial
canonical class $\KX$ and trivial Euler class $e$, the ideal
$\mathfrak B_X$ is the entire cohomology ring $H^*(X)$, and thus
the algebra $\Wbx$ coincides with $\Wax$. In addition,
Theorem~\ref{th_main} implies that the field $\jj^p(\alpha)(z)
=\sum_{n \in\Z} \jj^p_n(\alpha) z^{-n-p-1}$ is a {\rm primary}
field of conformal weight $(p+1)$ with respect to the Virasoro
field.

(ii) For a general projective surface $X$, we note that the
leading term in the commutators of the algebra $\Wax$ given in
Theorem~\ref{th_main} is precisely the $\W$ algebra introduced in
Sect.~\ref{sect_walg} associated to the ring $H^*(X)$. Therefore,
we are justified to regard the $\W$ algebras $\Wax$ and $\Wbx$ in
general as certain topological deformation of the $\W_{1+\infty}$
algebra in the framework of Hilbert schemes.
\end{remark}
\section{\bf Proof of Theorem~\ref{th_main}}
\label{sect_proof}

\subsection{Outline of the proof}  \label{outline}
$\,$

We denote the three terms in (\ref{eq_field}) by
${}^1\jj^p_m(\alpha), {}^2\jj^p_m(\alpha)$ and
${}^3\jj^p_m(\alpha)$ respectively, that is, $\jj^p_m(\alpha)
={}^1\jj^p_m(\alpha)+ {}^2\jj^p_m(\alpha)+ {}^3\jj^p_m(\alpha)$.
Note that $[{}^i \jj^p_m(\alpha),{}^j \jj^q_n(\beta)] =0$ $(i, j
=2,3)$ by the transfer property and $e^2 =0$. Thus,
\begin{eqnarray}  \label{eq_summand}
[\jj^p_m(\alpha), \jj^q_n(\beta)] \nonumber
 &=& [{}^1\jj^p_m(\alpha), {}^1\jj^q_n(\beta)]
     +[{}^1\jj^p_m(\alpha), {}^2\jj^q_n(\beta)]
     +[{}^1\jj^p_m(\alpha), {}^3\jj^q_n(\beta)]\nonumber \\
 & & + [{}^2\jj^p_m(\alpha), {}^1\jj^q_n(\beta)]+
     [{}^3\jj^p_m(\alpha), {}^1\jj^q_n(\beta)].
\end{eqnarray}

Recall that the commutation relations of Heisenberg generators can
be recasted equivalently in terms of the operator product
expansion (OPE) as (cf. \cite{Kac}; \cite{FB}, section 3.3):
\begin{eqnarray}  \label{eq_boson}
a(\alpha)(z)\; a(\beta)(w) \sim \frac{-\int_X(\alpha
\beta)}{(z-w)^2}.
\end{eqnarray}
Here and below $\sim$ means that the regular terms with
respect to $(z-w)$
on the right hand side of the OPEs are omitted; it is well
known the meromorphic terms carry all the information about the
corresponding commutation relations.

We derive from (\ref{eq_boson}) the following two OPEs:
\begin{eqnarray}   \label{eq_taylor}
a(\alpha)(z)\; \pa a(\beta)(w) \sim \frac{-2\int_X(\alpha
\beta)}{(z-w)^3}, \quad a(\alpha)(z)\; \pa^2 a(\beta)(w) \sim
\frac{-6\int_X(\alpha \beta)}{(z-w)^4}.
\end{eqnarray}

Our goal in this section is to compute the commutator
$[\jj^p_m(\alpha), \jj^q_n(\beta)]$ by using the OPE technique. To
this end, we need to compute the following OPEs

\begin{eqnarray}
 &(:a(z)^{p+1}:(\ta)) (:a(w)^{q+1}:(\tb)),& \label{eq_opeone} \\
 &(:a(z)^{p+1}:(\ta)) (:a(w)^{q-1}:(\teb)),& \label{eq_opetwo}\\
 &(:a(z)^{p+1}:(\ta)) (:\pa^2a(w)a(w)^{q-2}:(\teb)),& \label{eq_opethree}
\end{eqnarray}
which correspond to the first three commutators on the
right-hand-side of (\ref{eq_summand}) respectively. Then we will
be able to rewrite the terms appearing in these OPEs as
combinations of (the derivative fields of) the four fields
$:a(z)^{p+q}:(\tab)$, $:a(z)^{p+q-2}:(\teab)$,
$:\pa^2(a(z))a(z)^{p+q-3}:(\teab)$, and
$:\pa^3(a(z))a(z)^{p+q-3}:(\teab)$, by using the following simple
lemma.

\begin{lemma}  \label{lem_elem}
As in the previous section, write the fields $a^r(z), \pa^s a(z)$,
etc simply as $a^r, \pa^s a$, etc respectively. Then for $N \ge
0$, we have
\begin{enumerate}
\item[{\rm (i)}] $\pa^2 :a^N: \,\,\,= \,\,\, N:(\pa^2 a) a^{N-1}:
  + N(N-1):(\pa a)^2 a^{N-2}:$

\item[{\rm (ii)}] $\pa^3 :a^N: \,\,\,= \,\,\, N:(\pa^3 a) a^{N-1}:
  +6{N \choose 2} :(\pa a) (\pa^2 a) a^{N-2}:
  + 6{N \choose 3} :(\pa a)^3 a^{N-3}:$

\item[{\rm (iii)}] $\pa :a^N \pa^2 a:
  \,\,\,= \,\,\, :(\pa^3a) a^N: + N:(\pa a) (\pa^2 a) a^{N-1}:$

\item[{\rm (iv)}] $\pa :(\pa a)^2 a^N:
     \,\,\,= \,\,\, 2 :(\pa a) (\pa^2 a) a^N: +N :(\pa a)^3 a^{N-1}:$

\item[{\rm (v)}] $\pa^3 :a^N: \,\,\,= \,\,\, 3N \pa:(\pa^2 a) a^{N-1}:
   -2N :(\pa^3 a) a^{N-1}: +6{N \choose 3} :(\pa a)^3 a^{N-3}:$.
\end{enumerate}
\end{lemma}
\begin{demo}{Proof}
Using the Leibnitz rule, we obtain
\begin{eqnarray*}
\pa^2 :a^N: = \pa (N :(\pa a) a^{N-1}:) = N :(\pa^2 a) a^{N-1}:
+N(N-1) :(\pa a)^2 a^{N-2}:.
\end{eqnarray*}
This proves part (i). The proofs of the remaining formulas are
similar.
\end{demo}

We will calculate various OPEs by using the Wick formula. Note
that the Wick formula in the vertex algebra literature (cf. 
\cite{Kac}; \cite{FB},  Lemma~11.2.6, pp. 193) 
has to be modified in order to apply to our setup.
First of all, we naturally incorporate the transfer property,
Lemma~\ref{tau_k_tau_{k-1}}~(i). Next, the contraction of two
pairs of free fields will give rise to a term depending on
$e\alpha\beta$ (cf. Lemma~\ref{tau_k_tau_{k-1}}~(iii)). So the
contraction of more than two pairs of fields will be zero thanks
to $e^2=0$ and the transfer property. We refer to 
\cite{FB}, pp. 193, for the terminology ``contraction'' used above.

Now we assume $p+q >4$ in subsections \ref{proof_eq_opeone},
\ref{proof_eq_opetwo}, \ref{proof_eq_opethree} and
\ref{bigger_than_4}.
\subsection{Calculation of the OPE (\ref{eq_opeone})}
\label{proof_eq_opeone} $\,$

We calculate the OPE (\ref{eq_opeone}) by using the Wick theorem
and (\ref{eq_boson}) as follows:
\begin{eqnarray*}
 &&    \left(:a(z)^{p+1}:(\ta)\right) \left( :a(w)^{q+1}:(\tb) \right) \nonumber \\
 &\sim& (p+1)(q+1):a(z)^pa(w)^q:(\tab)\frac{-1}{(z-w)^2} \nonumber \\
 &&     +2{p+1 \choose 2}{q+1 \choose 2} :a(z)^{p-1}
           a(w)^{q-1}:(\teab)\frac{1}{(z-w)^4} \nonumber
\end{eqnarray*}
which by the Taylor expansion at $(z-w)$ is
\begin{eqnarray}
 &\sim& -(p+1)(q+1):a(w)^{p+q}:(\tab)\frac{1}{(z-w)^2} \nonumber \\
 &&  -\frac{(p+1)p(q+1)}{p+q} \pa \left( :a(w)^{p+q}:\right)(\tab)\frac{1}{z-w}
 \nonumber \\
 && +2{p+1 \choose 2}{q+1 \choose 2} :a(w)^{p+q-2}:(\teab)\frac{1}{(z-w)^4} \nonumber \\
 && +\frac{(p+1)p(p-1)(q+1)q}{2(p+q-2)} \pa
 \left(:a(w)^{p+q-2}:\right)(\teab)\frac{1}{(z-w)^3}\nonumber \\
 && + {p+1 \choose 3}{q+1 \choose 2} \left [ 3 :\pa
   \left(:\pa a(w)\; a(w)^{p-2}:\right)a(w)^{q-1}:(\teab)\frac{1}{(z-w)^2}
   \right .\nonumber \\
 && + \left . :\pa^2
 \left(:\pa a(w)\; a(w)^{p-2}:\right)a(w)^{q-1}:(\teab)\frac{1}{z-w} \right ].
 \label{eq_ope1}
\end{eqnarray}
Applying Lemma~\ref{lem_elem} to the two terms in the square
brackets in (\ref{eq_ope1}), we obtain
\begin{eqnarray}
 && :\pa \left(:\pa a(w)\; a(w)^{p-2}:\right)a(w)^{q-1}:(\teab) \nonumber \\
 &=& \frac{p-2}{(p+q-2)(p+q-3)} \pa^2 \left(:a(w)^{p+q-2}:\right)(\teab) \nonumber \\
 &&  + \frac{q-1}{p+q-3} :\pa^2 a(w)\; a(w)^{p+q-3}: (\teab) \label{eq_ope11}
 \end{eqnarray}
and
\begin{eqnarray}
 && :\pa^2 \left(:\pa a(w)\; a^{p-2}(w):\right) a(w)^{q-1}:(\teab) \nonumber \\
 &=& \frac{(p-2)(p-3)}{(p+q-2)(p+q-3)(p+q-4)} \pa^3(:a(w)^{p+q-2}:)(\teab) \nonumber \\
 && + \frac{3(p-2)(q-1)}{(p+q-3)(p+q-4)} \pa(:\pa^2 a(w)\; a(w)^{p+q-3}:) (\teab)\nonumber \\
 && -\frac{(p-q)(q-1)}{(p+q-3)(p+q-4)} :\pa^3a(w)\; a(w)^{p+q-3}: (\teab).
 \label{eq_ope12}
\end{eqnarray}
\subsection{Calculation of the OPE (\ref{eq_opetwo})}
\label{proof_eq_opetwo} $\,$

We calculate the OPE (\ref{eq_opetwo}) by using the Wick theorem,
(\ref{eq_boson}) and Lemma~\ref{lem_elem}:

\begin{eqnarray*}
 && \left(:a(z)^{p+1}:(\ta)\right) \left(:a(w)^{q-1}:(\teb)\right)\nonumber  \\
 &\sim& (p+1)(q-1) :a(z)^p a(w)^{q-2}:(\teab) \frac{-1}{(z-w)^2}
\end{eqnarray*}
which by the Taylor expansion at $(z-w)$ is
\begin{eqnarray}
 &\sim& -(p+1)(q-1) :a(w)^{p+q-2}:(\teab) \frac{1}{(z-w)^2}\nonumber  \\
 && -\frac{(p+1)p(q-1)}{p+q-2} \pa \left(:a(w)^{p+q-2}: \right) (\teab) \frac{1}{z-w} .
 \label{eq_ope2}
\end{eqnarray}
\subsection{Calculation of the OPE (\ref{eq_opethree})}
\label{proof_eq_opethree} $\,$

We calculate the OPE (\ref{eq_opethree}) by using the Wick
theorem, (\ref{eq_boson}) and (\ref{eq_taylor}) as follows:

\begin{eqnarray*}
 && \left(:a(z)^{p+1}:(\ta)\right) \left(:\pa^2a(w)a(w)^{q-2}:(\teb)\right)\nonumber \\
 &\sim& (p+1)(q-2) :a(z)^p\; \pa^2a(w)\, a(w)^{q-3}: (\teab)
        \frac{-1}{(z-w)^2} \nonumber \\
 &&  +(p+1) :a(z)^p a(w)^{q-2}: \frac{-6}{(z-w)^4} \end{eqnarray*}
which by the Taylor expansion at $(z-w)$ is
\begin{eqnarray}
 &\sim& -(p+1)(q-2) : \pa^2a(w) a(w)^{p+q-3}: (\teab)
        \frac{1}{(z-w)^2} \nonumber \\
 && -(p+1)p(q-2) :\pa a(w)\, \pa^2a(w)\, a(w)^{p+q-4}: (\teab)
        \frac{1}{z-w} \nonumber \\
 && -6(p+1) :a(w)^{p+q-2}:(\teab)\frac{1}{(z-w)^4} \nonumber \\
 && -6(p+1)p :\pa a(w)\, a(w)^{p+q-3}:(\teab)\frac{1}{(z-w)^3} \nonumber \\
 && -3(p+1)p :\pa \left(\pa a(w)\, a(w)^{p-1}\right) a(w)^{q-2}:(\teab)\frac{1}{(z-w)^2} \nonumber \\
 && - (p+1)p :\pa^2 \left( \pa a(w)\,a(w)^{p-1}\right)\, a(w)^{q-2}:(\teab)\frac{1}{z-w} .
 \label{eq_ope3}
\end{eqnarray}

We denote by $B_i$ $ ( i =1, 2)$ the coefficients of
$\frac1{(z-w)^i}$ on the right-hand-side of the above equation
(\ref{eq_ope3}). We then rewrite them by using
Lemma~\ref{lem_elem} as follows:

\begin{eqnarray}
B_2
 &=& -(p+1)(3p+q-2) :\pa^2 a(w) a(w)^{p+q-3}:(\teab) \nonumber \\
  &&   -3(p+1)p(p-1) :(\pa a(w))^2a(w)^{p+q-4}: (\teab)\nonumber \\
 &=& - \frac{3(p+1)p(p-1)}{(p+q-2)(p+q-3)} \pa^2 \left(:a(w)^{p+q-2}:\right) (\teab) \nonumber \\
 &&  - \frac{(p+1)(4p+q-3)(q-2)}{(p+q-3)} : \pa^2 a(w)\, a(w)^{p+q-3}: (\teab).
 \label{eq_ope31}
\end{eqnarray}
and
\begin{eqnarray}
B_1
 &=& -(p+1)p
    \left[ (3p+q-5) :\pa a(w)\, \pa^2 a(w)\, a(w)^{p+q-4}: (\teab) \right.  \nonumber \\
 && \qquad       + :\pa^3 a(w)\, a(w)^{p+q-3}:(\teab) \nonumber \\
 && \qquad  + \left. (p-1)(p-2) :(\pa a(w))^3 a(w)^{p+q-5}: (\teab) \right] \nonumber \\
 &=& -(p+1)p \cdot
    \left[ {(p-1)(p-2) \over \prod_{i=2}^4 (p+q-i)} \pa^3
    \left(:a(w)^{p+q-2}:\right) (\teab)\right. \nonumber \\
 && \qquad +\frac{(q-2)(4p+q-7)}{\prod_{i=3}^4 (p+q-i)} \pa
    \left(:\pa^2 a(w)\, a(w)^{p+q-3}:\right)(\teab)\nonumber \\
 && \qquad -\left. \frac{2(p-1)(p-2)}{\prod_{i=3}^4 (p+q-i)}
    :\pa^3 a(w)\, a(w)^{p+q-3}: (\teab)  \right].  \label{eq_ope32}
\end{eqnarray}
\subsection{Calculation of the commutator $[\jj^p_m(\alpha), \jj^q_n(\beta)]$
when $p+q > 4$} \label{bigger_than_4} $\,$

It follows from (\ref{eq_summand}) and Lemma~\ref{lem_fourier}
that the commutator $[\jj^p_m(\alpha), \jj^q_n(\beta)]$ is given
by the coefficient of $z^{-m-p-1}w^{-n-q-1}$ of the following OPE:
\begin{eqnarray*}
 && \frac1{(p+1)(q+1)} (:a(z)^{p+1}:(\ta)) (:a(w)^{q+1}:(\tb)) \\
&+& \frac{q}{24 (p+1)} (2q+3n-n^2) (:a(z)^{p+1}:(\ta))(:a(w)^{q-1}:(\teb)) w^{-2}\\
&-& \frac{q(q-1)}{24 (p+1)} (:a(z)^{p+1}:(\ta))(:\pa^2 a(w)a(w)^{q-2}:(\teb)) \\
&+& \frac{p}{24 (q+1)} (2p+3m-m^2) (:a(z)^{p-1}:(\tea))(:a(w)^{q+1}:(\tb)) z^{-2} \\
&-& \frac{p(p-1)}{24 (q+1)} (:\pa^2
a(z)a(z)^{p-2}:(\tea))(:a(w)^{q+1}:(\tb)).
\end{eqnarray*}
Note that the contributions of the last two OPEs in the above can
be obtained (up to a sign) from the preceding two OPEs by
switching $p,m,z$ with $q,n,w$ simultaneously.
The first three OPEs have been computed in subsections
\ref{proof_eq_opeone}, \ref{proof_eq_opetwo} and
\ref{proof_eq_opethree} respectively. We will need the following
simple fact: for $|z|>|w|$, $i \ge 0$,
\begin{eqnarray*}
\frac1{(z-w)^{i+1}} = \sum_{k\ge 0} {k \choose i} w^{k-i}
z^{-k-1}.
\end{eqnarray*}

We observe from the earlier computations of OPEs that there will
be four possible terms in the commutator $[\jj^p_m(\alpha),
\jj^q_n(\beta)]$, that is, $:a^{p+q}:_{m+n}(\tab)$,
$:a^{p+q-2}:_{m+n}(\teab)$, $:(\pa^2 a)\,a^{p+q-3}:_{m+n}(\teab)$,
and $:(\pa^3 a)\,a^{p+q-3}:_{m+n}(\teab)$. We can now determine
the coefficients of these terms one by one from the computations
of OPEs (cf. (\ref{eq_ope1}), (\ref{eq_ope11}), (\ref{eq_ope12}),
(\ref{eq_ope2}), (\ref{eq_ope3}), (\ref{eq_ope31}),
(\ref{eq_ope32})).

The coefficient of $:a^{p+q}:_{m+n}(\tab)$ in the commutator
$[\jj^p_m(\alpha), \jj^q_n(\beta)]$ is:
\begin{eqnarray*}
\frac1{(p+1)(q+1)}
  \left [ -(p+1)(q+1)(m+p) -\frac{(p+1)p(q+1)}{p+q} (-m-n-p-q)
  \right ]
\end{eqnarray*}
where the two terms in the brackets come from the first two terms
in (\ref{eq_ope1}) respectively. So the coefficient of
$:a^{p+q}:_{m+n}(\tab)$ in $[\jj^p_m(\alpha), \jj^q_n(\beta)]$
equals
\begin{eqnarray} \label{eq_coeff}
\frac{pn-qm}{(p+q)}.
\end{eqnarray}

Similarly, the coefficient of $:(\pa^3 a)\, a^{p+q-3}:_{m+n}
(\teab)$ in $[\jj^p_m(\alpha), \jj^q_n(\beta)]$ is
\begin{eqnarray*}
&&-\frac1{(p+1)(q+1)} \frac{(p+1)p(p-1)(q+1)q}{12}
   \frac{(p-q)(q-1)}{(p+q-3)(p+q-4)} \\
&&- \frac{q(q-1)}{24(p+1)} (p+1)p \frac{2(p-1)(q-2)}{(p+q-3)(p+q-4)} \\
&&+ \frac{p(p-1)}{24(q+1)} (q+1)q
\frac{2(q-1)(p-2)}{(p+q-3)(p+q-4)}
\end{eqnarray*}
which is $0$. So $:(\pa^3 a)\, a^{p+q-3}:_{m+n} (\teab)$ does not
appear in $[\jj^p_m(\alpha), \jj^q_n(\beta)]$.

Finally, similar lengthy computations together with MAPLE show
that the coefficient of $:(\pa^2 a)\, a^{p+q-3}:_{m+n}(\teab)$ in
$[\jj^p_m(\alpha), \jj^q_n(\beta)]$ is
\begin{eqnarray} \label{eq_coeffsq}
 -\frac{(pn-qm)(p+q-1)(p+q-2)}{24},
\end{eqnarray}
and the coefficient of $:a^{p+q-2}:_{m+n}(\teab)$ in
$[\jj^p_m(\alpha), \jj^q_n(\beta)]$ is equal to
\begin{eqnarray} \label{eq_coeffok}
\Theta_{m,n}^{p,q} +\frac{\Omega_{m,n}^{p,q}}{12(p+q-2)}.
\end{eqnarray}
In the above, $\Omega_{m,n}^{p,q}$ is given in (\ref{eq_number}),
and $\Theta_{m,n}^{p,q}$ is defined to be
\begin{eqnarray*}
\frac{(pn-qm)(p+q-1)}{24} (2p+2q-2+3m+3n-(m+n)^2).
\end{eqnarray*}

Collecting all the above computations together
(Eqs.~(\ref{eq_coeff}), (\ref{eq_coeffsq}), and
(\ref{eq_coeffok})) and using Definition~\ref{J}~(i), we have
established Theorem~\ref{th_main} for $p+q>4$.

\subsection{Completion of the proof of Theorem~\ref{th_main}}
\label{completion} $\,$

For the cases of $(p,q)$ with $p+q \le 4$, we can proceed either
directly or using OPEs as before while keeping in mind that many
of the formulas can be simplified. That is, some of the equations
such as (\ref{eq_ope11}), (\ref{eq_ope12}), (\ref{eq_ope31}),
(\ref{eq_ope32}) are no longer needed, and in fact they do not
make sense as some of the denominators would be zero. This is the
reason why we need to treat these cases separately.
When the unordered pair $(p,q)$ satisfies $3 \le p+q \le 4$, we
have checked that the commutators are given by the same formula as
the one for $p+q>4$. Finally, when $0 \le p+q \le 2$, the
commutators are found to be those stated in the
Theorem~\ref{th_main}. (in fact, the formulas for the unordered
pairs $(0,0)$, $(1, 0)$, $(1, 1)$ are precisely Theorem
\ref{commutator}~(i), (ii), (iv) respectively). We omit the
tedious details here.


\begin{thebibliography}{ABCD}



\bibitem[Bor]{Bor} R.~Borcherds,
{\em Vertex algebras, Kac-Moody algebras, and the Monster}, Proc.
Natl. Acad. Sci. USA {\bf 83} (1986) 3068--3071.



\bibitem[FB]{FB} E.~Frenkel and D.~Ben-Zvi,
{\em Vertex algebras and algebraic curves}, Math. Surveys and
Monographs {\bf 88}, Amer. Math. Soc. (2001).



\bibitem[FKRW]{FKRW} E.~Frenkel, V.~Kac, A. Radul and W. Wang,
{\em ${\mathcal W}_{1+\infty}$ and ${\mathcal W}(gl_N)$ with
central charge~$N$}, Commun. Math. Phys. {\bf 170} (1995)
337--357.



\bibitem[FW]{FW} I. Frenkel and W. Wang,
{\em Virasoro algebra and wreath product convolution}, J.~Alg.
{\bf 242} (2001) 656-671.



\bibitem[Got]{Got} L. G\"ottsche,
{\em The Betti numbers of the Hilbert scheme of points on a smooth
projective surface}, Math. Ann. {\bf 286} (1990) 193--207.



\bibitem[Gro]{Gro} I.~Grojnowski,
{\em Instantons and affine algebras I: the Hilbert scheme and
vertex operators}, Math. Res. Lett. {\bf 3} (1996) 275--291.



\bibitem[Kac]{Kac} V. Kac,
{\em Vertex Algebras for Beginners}, Second Edition, Univ. Lect.
Ser. {\bf 10}, Amer. Math. Soc. (1998).



\bibitem[LT]{LT} A. Lascoux and J.-Y. Thibon,
{\em Vertex operators and the class algebras of symmetric groups},
Preprint, math.CO/0102041.



\bibitem[Leh]{Leh} M. Lehn,
{\em Chern classes of tautological sheaves on Hilbert schemes of
points on surfaces}, Invent. Math. {\bf 136} (1999) 157--207.



\bibitem[LS1]{LS1} M. Lehn and C. Sorger, {\em Symmetric groups and
the cup product on the cohomology of Hilbert schemes}, Duke Math.
J. {\bf 110} (2001), 345--357.


\bibitem[LS2]{LS2} M. Lehn and C. Sorger,
{\em The cup product of the Hilbert scheme for $K3$ surfaces},
Preprint, math.AG/0012166.



\bibitem[LQW1]{LQW1} W.-P. Li, Z. Qin and W. Wang, {\em Vertex algebras and the
cohomology ring structure of Hilbert schemes of points on
surfaces}, Math. Ann. (to appear), math.AG/0009132.



\bibitem[LQW2]{LQW2} W.-P. Li, Z. Qin and W. Wang, {\em Generators
for the cohomology ring of Hilbert schemes of points on surfaces},
Intern. Math. Res. Notices No. {\bf 20} (2001) 1057--1074,
math.AG/0009167.



\bibitem[LQW3]{LQW3} W.-P. Li, Z. Qin and W. Wang, {\em Universality and
stability of cohomology rings of Hilbert schemes of points on
surfaces}, Preprint, math.AG/0107139.



\bibitem[Na1]{Na1} H. Nakajima,
{\em Heisenberg algebra and Hilbert schemes of points on
projective surfaces}, Ann. Math. {\bf 145} (1997) 379--388.



\bibitem[Na2]{Na2} H. Nakajima,
{\em Lectures on Hilbert schemes of points on surfaces}, Univ.
Lect. Ser. {\bf 18}, Amer. Math. Soc. (1999).


\bibitem[QW]{QW} Z. Qin and W. Wang, 
{\em Hilbert schemes and symmetric products:
a dictionary}, Proceedings for
the Conference ``Mathematical Aspects of Orbifold String Theory'',
Madison, Wisconsin, 2001, math.AG/0112070.


\bibitem[VW]{VW} C. Vafa and E. Witten,
{\em A strong coupling test of $S$-duality}, Nucl. Phys. {\bf B
431} (1994) 3--77.



\end{thebibliography}
\end{document}